\documentclass [12pt]{article}
\usepackage{amsmath}
\usepackage{amssymb}
\usepackage{amscd}
\usepackage{amsthm}
\usepackage{amsopn}
\usepackage{xspace}
\usepackage{verbatim}
\usepackage{amsmath}
\usepackage{amsfonts}
\newtheorem{theorem}{Theorem}[section]
\newtheorem{lemma}{Lemma}[section]
\newtheorem{proposition}{Proposition}[section]
\newtheorem{corollary}{Corollary}[section]
\newtheorem{definition}{Definition}[section]
\newtheorem{example}{Example}[section]
\newtheorem{acknowledgment*}{Acknowledgment}
\newtheorem{remark}{Remark}[section]
\numberwithin{equation}{section}
\newcommand{\rlemma}[1]{Lemma~\ref{#1}}
\newcommand{\rth}[1]{Theorem~\ref{#1}}

\newcommand{\ttt}{\mathbb{T}^n}

\newcommand{\vp}{\vec{p}}
\newcommand{\vdotx}{\dot{\vec{x}}}
\newcommand{\vx}{x}
\newcommand{\vy}{y}
\newcommand{\vz}{z}
\newcommand{\vlambda}{{\vec{P}}}
\newcommand{\slambda}{{P}}

\newcommand{\oM}{{\overline{\cal M}}}

\newcommand{\be}{\begin{equation}}
\newcommand{\ee}{\end{equation}}
\newcommand{\bd}{\begin{displaymath}}
\newcommand{\ed}{\end{displaymath}}

\newcommand{\oF}{F}
\newcommand{\oE}{\widehat{H}}

\newcommand{\vj}{{\vec{J}}}
\newcommand{\sJ}{{j}}
\newcommand{\dx}{dx}
\newcommand{\eps}{\varepsilon}

\newcommand{\R}{\mathbb R}
\newcommand{\M}{{\cal M}}
\newcommand{\Ss}{\mathbb{S}}

\newcommand{\CT}{C_0^1\left({\ttt}\times (0,T)\right)}

\renewcommand{\vec}[1]{\boldsymbol{#1}}

\begin{document}

\Large
\begin{center}{\bf Minimizers of Dirichlet functionals on the $n-$torus and the weak KAM Theory }\end{center}
\normalsize
\begin{center} G. Wolansky\end{center}

\begin{center}{\bf Resum\'e } \end{center}
Etant donn\'e  une mesure  $\mu$ sur le tore $n$-dimensionnel $\ttt$
et un vecteur de rotation  $\vlambda\in\R^n$, on \'etudie la
question de l'existence d'un minimiseur pour l'int\'egrale
$\int_{{\ttt}} |\nabla\phi+\vlambda|^2 d\mu$. Ce probl\`eme conduit
naturellement \`a une classe d'\'equations aux d\'eriv\'ees
partielles elliptiques et \`a une classe de probl\`emes de transport
optimal (Monge-Kantorovich) sur le tore. Il est aussi li\'e \`a la
th\'eorie d'Aubry-Mather en dimension sup\'erieure, qui traite les
ensembles invariants pour  des Lagrangiens p\'eriodiques, connue
sous le nom de  th\'eorie KAM faible.

\begin{center}{\bf Abstract} \end{center}
Given a probability measure $\mu$ on the $n-$torus $\ttt$ and a
rotation vector $\vlambda\in\R^n$, we ask whether  there exists a
minimizer to the integral $\int_{{\ttt}} |\nabla\phi+\vlambda|^2
d\mu$. This problem leads, naturally, to a class of elliptic PDE and
to an optimal transportation (Monge-Kantorovich) class of problems
on the torus. It  is also related to higher dimensional Aubry-Mather
theory, dealing with invariant sets of periodic Lagrangians, and is
known as the "weak-KAM theory".
\newpage
\tableofcontents \newpage

\section{Overview}
\subsection{Motivation}
Consider the functional \be\label{GAfunc}
H^\eps_{\vlambda}(u):= \frac{\eps^2}{2}\int_{{\ttt}} |\nabla u +
i\eps^{-1}\vlambda u|^2 dx- G(|u|) \ee where $\ttt:= \R^n\mod
\mathbb{Z}^n$ is the flat $n-$torus, $\vlambda\in\R^n$ a prescribed,
constant  vector, $u\in \mathbf{W}^{1,2}(\ttt; \mathbb{C})$ is
normalized via $\int_{{\ttt}}|u|^2 dx=1$ and $G$ is convex (possibly
non-local) functional of  $|u|$. A critical point $u$  of
$H^\eps_{\vlambda}$ can be considered as a periodic function on
$\R^n$. The function \be\label{bloch}u_0(x)= e^{i\vlambda\cdot
x/\eps} u(x)\ee is considered as a function on $\R^n$ as well.
 \vskip .2in\noindent {\bf Examples:}
\begin{description}
\item{i)} $G(|u|)= -\int_{{\ttt}} \Xi(x)|u|^2dx $ where $\Xi$ is a smooth
potential on $\ttt$. A critical point of (\ref{GAfunc}) is an
eigenvalues problem of the  Schr\"{o}dinger operator \\
$H_{\vlambda}:=- \eps^2(\nabla+i{\vlambda})^2  +2 \Xi $ on the
torus. The substitution (\ref{bloch})
  leads to a {\it Bloch state}  \be\label{el1} -\eps^2 \Delta u_0 + 2\Xi
u_0 + E u_0=0 \  \text{on} \ \R^n \ .  \ee
\item{ii)} \ Self-focusing  nonlinear  Schr$\ddot{\text{o}}$dinger equation. Here
$G(|u|)= \int_{{\ttt}}|u|^\sigma$ where $2< \sigma<2(2+n)/n$. An
extremum $u$ of $H^\eps_{\vlambda}$ with this choice is given by the
nonlinear eigenvalue problem for $u_0$ on $\R^n$: \be\label{el2}
\eps^2\Delta u_0 + \sigma|u_0|^{\sigma-1} u_0 = Eu_0
 \ . \ee
\item{iii)} \ The choice
$$G(|u|):=\sup_{V\in C^1(\ttt)}\left(
-\int_{{\ttt}}\frac{1}{2\gamma}|\nabla V|^2- V(|u|^2-1)\right) dx \
,
$$  leads  to the Schr$\ddot{\text{o}}$dinger Poisson system (see, e.g [BM], [AS], [IZL])
for  attractive (gravitational) field. Again, $u_0$  solves $$
-\eps^2\Delta u_0 - V u_0 =E_{(\vlambda)} u_0
$$ where $V$ is a periodic function on $\R^n$ solving $\Delta V +\gamma(|u_0|^2-1)=0
$.
\end{description}
 In addition to the (rather obvious) spectral asymptotic
questions, there are  additional motivations  for the study of
this problem, as described below. \\
 The short wavelength limit of
the reduced wave equation in a periodic lattice is described as
$$ \Delta u_0 +\eps^{-2} N(x)u_0=0 \ ,
$$ where $N(x)$ is the ($\eps$ independent) periodic function
representing the {\it refraction index} of the lattice and
$\eps\rightarrow 0$ stands for the wavelength. See, e.g. [RW].
Suppose one can measure the intensity $|u_0|$ and the direction
$\hat{\vlambda}:= \vlambda/|\vlambda|$ of the carrier wavenumber  of
an electromagnetic wave $u_0(x)=e^{i\vlambda \cdot x}u(x)$  in this
lattice (here, again, $u$ is periodic). Then $N$ can be recovered
from \be\label{Ndef} N\equiv
\frac{1}{2}|\nabla\phi+\hat{\vlambda}|^2 \ee where $\phi$ is the
minimizer of $F(\rho, \hat{\vlambda})$  for a normalized
$\rho=|u_0|^2$   (see (\ref{dirichlet}) below). Alternatively,
suppose we need to {\it design} a lattice for a prescribed
electromagnetic intensity $|u_0|^2$ and wave propagation
$\hat{\vlambda}$. Then (\ref{Ndef}) is the solution as well!

\subsection{The Effective Hamiltonian}

Note that, in  cases i-iii, I referred to {\it critical points} of
$H^\eps_{\vlambda}$. If one looks at minimizers of this functional,
or even critical points of finite ($\eps-$independent) Morse index,
then one may expect singular limits as $\eps\rightarrow 0$.
 However, there is a formal way to obtain  nonsingular
limits of these equations as $\eps\rightarrow 0$ as follows:\par
Substitute the WKB anzatz (see [K])  $u_\eps:=\sqrt{\rho}
e^{i\phi/\eps}$ in (\ref{GAfunc}), where $\phi\in C^1(\ttt)$ and
$\rho\in C^1(\ttt)$ is non-negative function satisfying
$\int_{{\ttt}}\rho =1$. Then  $$ \lim_{\eps\rightarrow 0}
H_{\vlambda}^\eps(u_\eps)= \frac{1}{2}\int_{{\ttt}}
|\nabla\phi+\vlambda|^2 \rho(x)dx - G(|\rho|^{1/2}) \ . $$ Now,
define \be\label{dirichlet}F(\rho,
\vlambda):=\frac{1}{2}\inf_{\phi\in C^1(\ttt)} \int_{{\ttt}}
|\nabla\phi+\vlambda|^2 \rho(x)dx\ee and \be \label{G0def}
H_G(\rho,\vlambda):= F(\rho, \vlambda) - G(|\rho|^{1/2}) \ . \ee
Note that $F$ is concave as a function of $\rho$ for fixed
$\vlambda\in\R^n$. Then we can (at least formally) look for the {\it
maximizer} of $H_G(\cdot,\vlambda)$ over the set of densities
$\rho\in \mathbb{L}^1(\ttt; \R^+)$. It is the {\it asymptotic energy
spectrum} $H$ associated with $H^\eps_{\vlambda}$:
\be\label{energy0} \oE_G(\vlambda):= \sup_\rho\left\{
H_G(\rho,\vlambda)  \ ; \ \ \rho\in \mathbb{L}^1(\ttt; \R^+) \ , \
\int_{{\ttt}}\rho=1\ . \right\} \ee
 The pair
$(\phi,\rho)$ which realizes the minimum (res. maximum) of $F$ (res.
$\oE_G(\cdot, \vlambda)$) corresponds to an asymptotic critical
point of $H^\eps_{\vlambda}$ and verifies the  Euler-Lagrange
equations: \be \label{ellipeq} \nabla\cdot\left[\rho\left(
\nabla\phi+\vlambda\right)\right]=0 \ , \ee\be\label{ellipeq1}
\frac{1}{2}|\nabla\phi+\vlambda|^2 - G_\rho= E \ee on $\ttt$, where
$G_\rho$ is the  Fr\`{e}chet derivative of $\rho\rightarrow
G(|\rho|^{1/2})$ and  $E$ is a Lagrange multiplier corresponding to
the constraint $\int_{{\ttt}}\rho=1$.
\par
Of particular interest  is  the linear case (\ref{el1}). Here
(\ref{G0def}, \ref{energy0}) are reduced into \be \label{G0Wdef}
H_\Xi(\rho,\vlambda):= F(\rho, \vlambda) + \int_{\ttt}\Xi\rho  \ .
\ee \be\label{energyXi} \oE_\Xi(\vlambda):= \sup_\rho\left\{
H_\Xi(\rho,\vlambda)  \ ; \ \ \rho\geq 0 \ , \ \int_{{\ttt}}\rho=1\
. \right\} \ee and (\ref{ellipeq1})  takes the form of the
Hamilton-Jacobi equation on the torus: \be\label{HJ}
\frac{1}{2}|\nabla\phi+\vlambda|^2 + \Xi= E \   \ \ee which is
independent of $\rho$ and so is decoupled from (\ref{ellipeq}).
\par
Suppose now there exists a maximizer  $\rho_0$  of (\ref{energyXi}).
Multiply (\ref{HJ}) by $\rho_0$ and integrate over $\ttt$ to obtain
\be\label{Gnomore}E=F(\rho_0,\vlambda) +\int_{{\ttt}}\Xi \rho_0 :=
\oE_\Xi( \vlambda)\  \ .  \ee
 In particular, the Lagrange multiplier $E$ is
identical to the asymptotic energy spectrum $\oE_\Xi(\vlambda)$.  An
important point to be noted, at this stage, is that the asymptotic
spectrum is in the oscillatory domain of the periodic
Schr$\ddot{\text{o}}$dinger equation, that is
$\oE_\Xi(\vlambda)\geq\max_{{\ttt}}\Xi$ necessarily holds for any
$\vlambda\in\R^n$, and $\oE_\Xi(\vlambda)>\max_{{\ttt}}\Xi$ if
$|\vlambda|$ is sufficiently large (see
Proposition~\ref{cordelta2}).
 \par
The function $\oE_\Xi=\oE_\Xi(\vlambda)$ defined in (\ref{Gnomore})
is considered by Evans and Gomes  [EG1, EG2] as the {\it Effective
Hamiltonian} corresponding to
$$h_\Xi(p,x)= |p|^2/2 + \Xi(x) \ . $$ If $\psi(x,\vlambda):=\phi(x)$
is a solution of the Hamilton-Jacobi equation (\ref{HJ})
corresponding to a given $\vlambda$, then a canonical change of
variables
$$ \vp=\nabla_x\psi +\vlambda \ \ ; \ \ \vec{X}=\nabla_{\vlambda}\psi+\vx$$
reduces the Hamiltonian equation to an integrable system defined by
the Hamiltonian $\oE_\Xi$, that is, $\vlambda$ is a cyclic variable
and hence a constant of motion.
\par
In general, such a solution does not exist for any
$\vlambda\in\R^n$. However, (\ref{energyXi}) suggests a way to
define the effective Hamiltonian $\oE_\Xi$ {\it without} the
assumption that (\ref{HJ}) is solvable. We note, at this stage, that
(\ref{energyXi}) seems to be the dual of \be\label{Hdual}
\oE_\Xi(\vlambda)= \inf_{\phi\in C^\infty(\ttt)}\sup_{x\in \ttt}
h_\Xi(\nabla\phi+\vlambda, x) \  \ee which was suggested by Gomes
and Oberman [GO] as a numerical tool for evaluating $\oE_\Xi$.
\subsection{Objectives}
As we shall see below,  the supremum in (\ref{energyXi}) is not
attained in $\mathbb{L}^1(\ttt)$, in general, but in the set $\oM$
of Borel probability measures on $\ttt$.  This, together with
(\ref{G0Wdef}),  motivates us to extend the domain  of $F$ in
(\ref{dirichlet}) from the set of non-negative densities in
$\mathbb{L}^1(\ttt)$  to $\oM$.  Similarly, the functional $H_\Xi$
(\ref{G0Wdef}) is extended to $\oM$ as well.  Our first object is
\vskip .2in\noindent {\bf I.} \ \ \ {\it to define a generalized
minimizer $\phi$ of $F$.} \vskip .2in
 The effective Hamiltonian
 plays a major rule  in
 the {\it weak KAM Theory}. See [F, GCJ, Man, Mat1, S] among other
 references.
For the convenience of the reader we review
 the fundamentals of the weak KAM Theory and Mather measures  in  section~\ref{secminmeas}.
 Our second object is
\vskip .2in\noindent {\bf II.} \ \  {\it to relate the functional
$H_\Xi$ to the weak KAM theory. In particular, to relate the
generalized minimizer $\phi$ of $F$ to the minimal Mather measure. }
\vskip .2in An excellent reference to the Monge-Kantorovich Theory
of optimal transportation is the book of Villani [V]. The relation
between M-K theory  and the weak KAM theory was
  suggested in [E] and  further elaborated in a series of
 publications, among which
[G,  LSG, BB]. Essentially, it relates the minimal (Mather) measures
of a given Lagrangian to a measure which minimizes a certain optimal
transportation plan. The third object of this paper is \vskip
.2in\noindent {\bf III.} \ \ {\it to approximate $F$ and $H_\Xi$
 by an optimal transportation
functions $F_T$  and $H_{\Xi,T}$ , respectively. } \vskip.2in
Finally, we use the suggested functionals to establish an
alternative to (\ref{Hdual}) for  the evaluation of
$\oE_\Xi(\vlambda)$: \vskip.2in\noindent {\bf IV.} {\it Establish a
combinatorial search algorithm for evaluation of $\oE_\Xi(\vlambda)$
to any degree of approximation.}

\section{Lagrangian dynamics on the torus}\label{secminmeas}
\subsection{The Aubrey-Mather Theory and minimal orbits}

  Let
\be\label{mechkag}L(\vp,\vx):= \frac{|\vp|^2}{2}-\Xi(\vx) \ , \ee a
Lagrangian function defined on $\R^n\times\R^n$ where  the potential
$\Xi$ which is $1-$ periodic in  all the variables $\vx=(x_1,\ldots
x_n)$. For a given orbit $(\vx(t),\vp(t))$ of the associated
Euler-Lagrange equation \be\label{Lageq} \dot{\vx}=\vp \ \ \ ; \ \ \
\dot{\vp}+\nabla_{\vx}\Xi=0  \ , \ee
 a {\it rotation vector} $\vj\in \R^n$ is
assigned to this orbit provided the limit
\be\label{classrot}\vj=\lim_{|t|\rightarrow\infty}t^{-1}\vx(t) \ \ee
exists.
 As a trivial  example, consider the $\vx-$independent Lagrangian
 where $\Xi\equiv 0$.
 Since $\vp$ is a constant of motion
and $\dot{\vx}=\vp$, the rotation vector is defined for each orbit
via $\vj=\vp$.  For general $\Xi$ the rotation vector is {\it not}
defined for any orbit, in general.  The object of the {\it
classical} KAM theory  (see, e.g. [GCJ]) is the study of  small
perturbation of an integrable system, e.g. for Lagrangians of the
form (\ref{mechkag}) where the potential $\Xi$ is small. In
particular, it studies  families of solutions of such systems which
preserve the rotation vector.
\par
In the eighties, Aubry [A] and Mather [Mat] (see also [Mo])
discovered that Lagrangian flows which induce a monotone, symplectic
twist maps on a two dimensional annulus, possess orbits of any given
rotation number (in the twist interval), even if the corresponding
Lagrangian is not  close to an integrable one. The characterization
of these orbits is variational: They are minimizers of the
Lagrangian action with respect to any local variation of the orbit.
In general, they are embedded in invariant tori of the Lagrangian
flow.
\par There is still another approach to
invariant tori of Lagrangian/Hamiltonian systems.  An invariant
Lagrangian torus can be obtained as a solution of the corresponding
Hamilton-Jacobi equation as follows: Suppose there exists, for some
$\vlambda\in\R^n$,  solution $\phi\in C^{1,1}(\R^n)$ which is
1-periodic in each of the coordinates $x_j$ of $\vx=(x_1,\ldots
x_n)$,  for
$$\frac{1}{2}|\nabla\phi(\vx)+\vlambda|^2+\Xi(\vx)=E \ \ ,  \ E\in \R \ . $$
Then the graph of the function $(\vx, \vlambda+\nabla\phi(\vx))$
represents an {\it invariant torus} of the Lagrangian flow
associated with $L$ [Man].  The projection on $\vx$ of any orbit in
this invariant set is obtained by a solution of the system
\be\label{dotx=nab}\dot{\vx}=\vlambda+\nabla\phi(\vx) \ .  \ee In
the case $n=2$ the rotation vector $\vj\in\R^2$ is defined as in
(\ref{classrot}) for any such orbit, given by (\ref{dotx=nab}).
\subsection{Weak KAM and minimal invariant measures}
  For dimension higher than 2,
there are counter-examples: There exists a Lagrangian system on the
3 dimensional torus, induced by a metric on this torus, for which
there are no minimal geodesics, save for a finite number of rotation
vectors [H]. Moreover, it is not known that the limit
(\ref{classrot}) exists for any orbit of  (\ref{dotx=nab}), if
$n>2$. Hence, an
 extension of Aubry-Mather theory  to higher
dimensions is not a direct one. If, however, we replace the notion
of an orbit  by an invariant measure, then it is  possible to extend
the Aubry-Mather theory to higher dimensions.   The relaxation of
orbits to invariant measures (and the corresponding minimal orbits
to minimal invariant measures) leads to
 the "weak KAM Theory".
 \par Let ${\cal M}_L$ be the set of all  probability measures on the tangent bundle $\ttt\times\R^n$  which
 are invariant with respect to the flow induced by the  Lagrangian $L$.  The {\it rotation
vector} $\vec{\alpha}:{\cal M}_L\rightarrow \R^n$  is
\be\label{rotnum} \vec{\alpha}(\nu):=\int_{\ttt\times\R^n} \vp
d\nu(x,\vp) \ , \ee and, for any $\alpha\in\R^n$,  the set of all
${\cal M}^{\vj}_L\subset {\cal M}_L$ corresponds to all $\nu\in{\cal
M}_L$ for which $\vec{\alpha}(\nu)=\vj$.
 \par
 A minimal measure
associated with a rotation vector $\alpha$ is defined by
\be\label{minmeas} \nu_{\vj}=\arg\min_{\nu\in {\cal M}^{\vj}_L}
\int_{\ttt\times\R^n}L d\nu\in {\cal M}_L^{\vj} \ . \ee Its dual
representation is given by minimizing the Lagrangian $L^{\vlambda}:=
L(\vp,x)-\vlambda\cdot \vp$ over the {\it whole set} of invariant
measures\footnote{Note that $\nu$ is an invariant measure of $L$ if
and only if it is an invariant measure of $L_{\vlambda}$ \ . }
${\cal M}_L$. The measure $\nu_\vlambda\in {\cal M}_L$ is called a
{\it Mather measure} if \be\label{minmeasL}
\nu_{\vlambda}=\arg\min_{\nu\in {\cal M}_L}
\int_{\ttt\times\R^n}L^{\vlambda} d\nu\in {\cal M}_L \ . \ee

These minimal measures are relaxations of minimal invariant orbits
of the Aubry-Mather theory. Their properties and the geometry of
their supports are the fundamental ingredients of the developing
weak KAM theory. For further details, see [Mat1], [Man], [F] or
consult
 [S] for applications and further references.
 \par
 It should
be stressed, however, that the investigation of the functional $\oF$
(\ref{dirichletmu}) carried in  the present paper is not restricted
to minimal (Mather) measures. In fact, Mather measures (and their
$\ttt$ projections) are defined only for smooth enough Lagrangian
systems which allow the existence of dynamics, e.g (\ref{mechkag})
where $\Xi\in C^{1,1}(\ttt)$. Since we are motivated, between other
things,  by quantum dynamics and wave equation, we must assume much
less, e.g the  Schr$\ddot{\text{o}}$dinger equation (\ref{el1}) is
well posed if the potential $\Xi$ is only continuous.

\section{An overview of the main results}
\subsection { List of
symbols and definitions}\label{list}
\begin{enumerate}
 \item \
$\ttt=\R^n/ \mathbb{Z}^n$ the $n-$dimensional flat torus. It is
parameterized by \\ $x=(x_1,\ldots x_n)\mod \mathbb{Z}^n$. The
Euclidian  distance $\| x-y\|_{{\ttt}}$ on $\ttt$ is defined as
$\min_{z\in \mathbb{Z}^n}|x-y-z|$, where $x,y\in\R^n$ and $|\cdot |$
is the Euclidian norm on $\R^n$.
\item  \ $
\ttt\times\R^n$ is the tangent bundle of $\ttt$, that is, $
\ttt\times\R^n:=\R^n\times \ttt$. A vector in $\R^n$  is denoted by
a bold letter, e.g. $\vec{v}$. The same symbol will also define a
vector field, that is, a section in $ \ttt\times\R^n$, e.g
$\vec{v}=\vec{v}(x)$.
\item  \ ${\cal M}(D)$ stands for the set of
all probability normalized Borel measures  $\mu$ on some metric
space $D$, subjected to the dual topology of $C(D)$: $|\mu|:=
\sup_{\phi\in C(D), |\phi|_\infty=1} \int_D \phi d\mu$.  We denote
$\oM:= {\cal M}(\ttt)$ as the set of all such measures on the torus
$\ttt$.
\item \  A  Borel map ${\bf S}:D_1\rightarrow D_2$ induces a
map ${\bf S}_\#: {\cal M}(D_1)\rightarrow {\cal M}(D_2)$, as
follows:
\\
 ${\bf S}_\#\mu (A):=\mu\left({\bf S}^{-1}(A)\right)$ for any Borel set $A\in D_2$.
${\bf S}_\#\mu$ is called the {\it push-forward} of $\mu\in {\cal
M}(D_1)$ into ${\cal M}(D_2)$.
\item $\pi:\mathbb{T}^n\times \R^n\rightarrow \mathbb{T}^n$ is the projection (natural
embedding) of $\mathbb{T}^n$ in $\mathbb{T}^n\times\R^n$, namely
$\pi(x,\vp)=x$ for $(x,\vp)\in \mathbb{T}^n\times\R^n$. In
particular, $\pi_\#:{\cal M}(\mathbb{T}^n\times \R^n)\rightarrow\oM$
\   so $\mu=\pi_\#\nu\in \oM$ is the $\mathbb{T}^n$ marginal of
$\nu\in {\cal M}(\mathbb{T}^n\times\R^n)$.
\item\label{orbits} $\oM_T:= C([0,T], \oM)$.
An element  $\hat{\mu}\in\oM_T$ is denoted by $\hat{\mu}:=
\mu_{(t)}$, $0\leq t\leq T$. For any $\mu_1, \mu_2 \in\oM$, the set
$\oM_T(\mu_1, \mu_2)\subset\oM$ is defined as all
$\hat{\mu}\in\oM_T$ for which $\mu_{(0)}=\mu_1, \mu_{(T)}=\mu_2$. If
$\mu_1=\mu_2\equiv \mu$ we denote $\oM_T(\mu):= \oM_T(\mu, \mu)$.
  \item\label{listm^2}  \ $\overline{\cal M}^{(2)}:= {\cal M}(\ttt\times
\ttt)$ and $\Pi_i:\ttt\times\ttt\rightarrow \ttt$ the projection  on
the $i$ factor, $i=1,2$. For $\mu_1\in\oM$, $\mu_2\in\oM$ define
$\overline{\cal M}^{(2)}(\mu_1, \mu_2)$ as the set of all $\sigma\in
\overline{\cal M}^{(2)}$ for which $\Pi_{i,\#}\sigma= \mu_i$,
$i=1,2$. $\overline{\cal M}^{(2)}(\mu):= \overline{\cal
M}^{(2)}(\mu, \mu)$.
 \item \  $C^1(\ttt)$ is the set of all $C^1$ smooth functions on
 $\ttt$.
 \item \  Recall the definition of a subgradient of a
function $h:B\rightarrow \R$ defined on a Banach space $B$:  For
$b\in B$, $$
\partial_b h:= \left\{ b^*\in B^* \ ,  h(b^{'})\geq h(b) + \langle b^{'}-b,
b^*\rangle \ \ \text{ For any } \ \ b^{'}\in B \ . \right\} $$
\end{enumerate}
\subsection{$\oF$ and its generalized minimizers}
There is a close relation between the minimal (Mather)  measures
described in Section~\ref{secminmeas} and the minimizer of the
function $F$ defined in (\ref{dirichlet}), where $\rho$ is the
density of the $\ttt$ marginal of a minimal measure.
 \par
 In general, however, there are no smooth densities to the marginals of minimal measures on $\ttt$.
 Motivated by this, we extend the definition of $F$
 to the set $\oM$ of {\it all} probability Borel measures on $\ttt$:
 \be\label{dirichletmu}F(\mu,
\vlambda):=\frac{1}{2}\inf_{\phi\in C^1(\ttt)} \int_{{\ttt}}
|\nabla\phi+\vlambda|^2 d\mu   \ee and   \be\label{legF}
\oF^*(\mu,\vj):= \sup_{\vlambda\in \R^n} \left[ \vlambda\cdot\vj -
\oF(\mu,\vlambda)\right] \  \ee its convex dual on $\R^n$.
\par
The first question we address is the existence of minimizers of
(\ref{dirichletmu}) for a general measure $\mu\in \oM$. Evidently,
there is no sense of solutions to the elliptic problem
(\ref{ellipeq}) for such $\mu$. Our first result, given in
\rth{th1},  indicates the existence and uniqueness of a minimizer in
a generalized sense (see Definition~\ref{defweak} below). In
\rth{th3} we discuss the relation between these minimizers and the
solutions of the elliptic problem (\ref{ellipeq}).
\subsection{The Effective Hamiltonian}
The second question concerns the  maximizers of (\ref{Gnomore}),
extended to the entire set $\oM$. Let \be \label{Gnomorebar}
\oE_\Xi(\vlambda):= \sup_{\mu\in\oM}\left[
\oF(\mu,\vlambda)+\int_{{\ttt}}\Xi d\mu \right] \ , \ee where
$\Xi\in C(\ttt)$.
 \rth{connect} relates the generalized minimizers of \rth{th1} to
 the minimal (Mather) measure associated
with the Lagrangian (\ref{mechkag}) corresponding to $\vlambda$ via
(\ref{minmeasL}). It claims that, if $\Xi\in C^{1,1}(\ttt)$, then
the generalized solution of $\oF(\mu,\vlambda)$ corresponding to
$\mu$ which maximizes (\ref{Gnomorebar}) is, indeed, a Mather
measure associated with the Lagrangian (\ref{mechkag}). In this
sense, the weak minimizers of \rth{th1} can be considered as
generalized Mather measures for Lagrangians with only continuous
potentials.
\subsection{On the continuity of $\oF$}
 The third question addressed is  the
 continuity
property of $\oF$ with respect to  $\mu$. It is rather easy to
observe that $\oF$ is convex in $\vlambda$ on $\R^n$, and concave in
$\mu$ on $\oM$. These imply that $\oF$ is continuous on $\R^n$, but
only {\it upper-semi-continuous} in the natural topology of $\oM$,
which is the weak-$*$ topology induced by $C^*(\ttt)$. That is, if
$\mu_j\rightarrow \mu$ in $C^*(\ttt)$, then \be\label{spappineq}
\lim_{j\rightarrow\infty} F( \mu_j, \vlambda)\leq F(\mu, \vlambda)\
\ \
 \ee
holds. In general, there is {\it no continuity} of $\oF$ over $\oM$
with the $C^*$ topology.
\par\noindent {\bf Examples:}
\begin{enumerate}
\item For any atomic measure $\mu=\sum m_i\delta_{x_i}\in \oM$, we
can easily verify that $\oF( \mu, \vlambda)\equiv 0$ for any
$\vlambda\in\R^n$. In particular, if $\mu_N$ is a sequence of {\it
empirical} measures: $\mu_N:= N^{-1}\sum_{i=1}^N \delta_{x_i}$
satisfying  $\mu_N\rightarrow \mu\in \oM$ in $C^*(\ttt)$, then the
inequality in (\ref{spappineq}) is strict whenever $\oF(\mu,
\vlambda)>0$.
\item
Let $n=1$, so $\ttt$ is reduced to the circle $\mathbb{S}^1$.
Suppose $\mu\in \oM( \mathbb{S}^1)$ admits a smooth density
$\mu(\dx)=\rho(x)dx$. If $\rho>0$ on $\Ss^1$ then the continuity
equation (\ref{ellipeq}) reduces to a constant
$\sJ=\rho(\phi_x+\slambda)$. This implies \be\label{facj} \oF(\mu,
\slambda)= \frac{1}{2}\int_{{\Ss}^1} \rho|\phi_x+\slambda|^2 dx =
\frac{\sJ^2}{2}\int_{{\Ss}^1} \rho^{-1}dx \  \ee as well as
$$ \int_{{\Ss}^1} \rho^{-1} dx = \int_{{\Ss}^1}\sJ^{-1}(\phi_x+\slambda)
dx = \frac{ \slambda}{j} \Longrightarrow \sJ=
\slambda\left(\int_{{\Ss}^1} \rho^{-1} dx\right)^{-1}  \ .   $$
Substitute in (\ref{facj}) to obtain
$$ \oF( \mu,\slambda)= \frac{|\slambda|^2}{2}\left(\int_{{\Ss}^1} \rho^{-1}
dx\right)^{-1} \ . $$

  In particular, $\int_{{\Ss}^1}\rho^{-1} dx=\infty$ iff $\oF(  \mu,\slambda)=0$ for $\slambda\not= 0$.   Any  sequence
$\mu_j(dx)=\rho_j(x)dx$ satisfying
$\int_{{\Ss}^1}\rho_j^{-1}=\infty$ which converges in $C^*(
{\Ss}^1)$ to $\mu(dx)=\rho(x)dx$ satisfying $\rho\in C^1({\Ss}^1)$,
$\rho>0$ on $\Ss^1$, is an example of strict inequality in
(\ref{spappineq}).
\end{enumerate}
\subsection{Lagrangian mappings}
In \rth{limtt} we show that $\oF$  can be approximated, as a
function on $\oM$,  by a {\it weakly continuous} function
$\oF_T(\cdot, \vlambda)$  which satisfies $\oF_T(\mu,
\vlambda)\rightarrow \oF(\mu,\vlambda)$  as $T\rightarrow 0$, for
any $\mu\in\oM$. For this, we represent an extension of $\oF$ to
{\it orbits} $\hat{\mu}: [0,T]\rightarrow \oM$ ,
$\left.\hat{\mu}\right|_{(t)}=\mu_{(t)}\in\oM$, $t\in [0,T]$, given
by
$$ F(\hat{\mu},\vlambda, T) :=\frac{1}{T}\inf_{\phi\in
\CT}\int_0^T\int_{{\ttt}}\left[\frac{\partial\phi}{\partial t}+
\frac{1}{2}|\nabla_x\phi+\vlambda|^2 \right]d\mu_{(t)} dt \ , $$ and
set $\oF_T(\mu,\vlambda)$ as the {\it supremum} of
$\oF(\hat{\mu},\vlambda,T)$ over all such orbits satisfying
$\mu_{(0)}=\mu_{(T)}=\mu$.  It is shown that $\oF_T(\mu, \vlambda)=
|{\vlambda}|^2/2 - D_{T\vlambda}(\mu)/(2 T^2)$ where
$D_{\vlambda}(\mu)$ is defined by the optimal Monge-Kantorovich
transport plant  from $\mu$ to itself, subjected to the cost
function $c(x,y):=\|x-y-\vlambda\|_{{\ttt}}^2$, where
$\|\cdot\|_{{\ttt}}$ is the Euclidian metric on $\ttt$. As an
example, consider the case $\mu=\delta_{x_0}$ for some $x_0\in
\ttt$. Then  $D_{\vlambda}(\delta_{x_0})= |\{\vlambda\}|^2$, where
$\{ \cdot\}$ stands for the fractional part $\{ \vlambda\} :=
\vlambda\mod \mathbb{Z}^n$. So
$$\oF_T(\delta_{x_0},\vlambda)= |\vlambda|^2/2-  |\{\vlambda T\}|
^2/(2T^2).    \ $$ If $T$ is sufficiently small so $\{ T\vlambda\}=
T\vlambda$ then $\oF_T(\delta_{x_0},\vlambda)=0$.

\subsection{Combinatorial search for the minimal measure}
 Our last object is to suggest an alternative to the
numerical algorithm for the calculation of the effective Hamiltonian
based on (\ref{Hdual}), introduced in [GO]. We take an advantage of
the following facts
\begin{description}
\item{i)}   An optimal transportation functional (such as $\oF_T(\mu,\vlambda)$) are continuous in the
weak topology of $\oM$.
\item{ii)}
The set of {\it empirical measures} is dense in the set of all
measures $\oM$.
\item{iii)} On the set of empirical measures of a fixed number of sampling points $j$, an optimal
transportation problem is reduced to a finite combinatorial problem
on the set of permutation on $\{1, \ldots j\}$ (Birkhoff's Theorem).
\end{description}
Applying  (i-iii) to the result of Theorem~\ref{limtt}, we
   obtain a discrete, combinatorial
algorithm for evaluating the effective Hamiltonian
$\oE_\Xi(\vlambda)$. This is summarized in Theorem~\ref{DisH}.

\section{Detailed description of the main results}\label{msec}
\subsection{Minimizers of the Dirichlet functional over the
$n-$torus}
 Let us recall the
definition, for $\vlambda\in\R^n$, $\vj\in\R^n$ and
$\mu\in\overline{\cal M}$: \be\label{Fbar2}
F(\mu,\vlambda):=\frac{1}{2}\inf_{\phi\in C^1(\ttt)}\int_{{\ttt}}
\left|\nabla\phi+\vlambda\right|^2 d\mu \ , \ee
 \be\label{F*bar2}
\oF^*(\mu,\vj):= \sup_{\vlambda\in \R^n} \left[ \vlambda\cdot\vj -
\oF(\mu,\vlambda)\right] \  \ee Let also ${\cal E}$ a function on
$\oM\times\R^n\times C^1(\ttt)$ defined as:
 \be\label{calE} {\cal E}(\mu,\vj,\phi):= \frac{1}{2}\left\{
\left| \vj-\int_{{\ttt}}\nabla\phi d\mu\right|^2-
\int_{{\ttt}}|\nabla\phi|^2d\mu \right\} \ . \ee

\par Next, we consider the notion of {\it weak solution} of
(\ref{ellipeq}), corresponding to the minimizer of (\ref{Fbar2}):
\begin{definition}:\label{deflift}
\begin{enumerate}
\item
The set $\Lambda \subset {\cal M}( \ttt\times \R^n)$ consists of all
probability measures $\nu(\dx d\vp)$ for which
$$\int_{ \ttt\times\R^n }\vp\cdot\nabla\theta(x) d\nu = 0 \ \ \
\forall\theta\in C^1(\ttt) \  \ \text{and} \  \int_{\ttt\times\R^n
}|\vp|^2d\nu<\infty \ .
$$
\item
Given  $\mu\in \oM$, the set $\Lambda_\mu\subset \Lambda$ of
liftings of $\mu$ is composed of all $\nu\in\Lambda$ for which
$\pi_\#\nu=\mu$.
\item For each $\vj\in\R^n$, the set $\Lambda_\mu^{\vj}$ is
defined as all  $\nu\in\Lambda_\mu$ which satisfies
$$ \int_{ \ttt\times\R^n} \vp d\nu = \vj \ . $$
\end{enumerate}
\end{definition}
\begin{remark}
The set $\Lambda_\mu$ is never empty. Indeed,
$\nu=\delta^{\vp}_0\otimes\mu\in\Lambda_\mu$ for any $\mu\in\oM$.
However, the set $\Lambda^{\vj}_\mu$ can be empty. For example, if
$\oF(\mu,\vlambda)=0$ for all $\vlambda\in\R^n$ then
$\Lambda^{\vj}_\mu=\emptyset$ for any $\vj\not= \vec{0}$.
\end{remark}
\begin{definition}\label{defweak}
For given $\vlambda\in\R^n$, $\nu\in\Lambda_\mu$ is
 a weak solution of
$\oF(\mu,\vlambda)$  provided
$$ \int_{ \ttt\times\R^n} \left[\frac{|\vp|^2}{2}-\vp\cdot\vlambda\right] d\nu \leq \int_{ \ttt\times\R^n}
 \left[\frac{|\vp|^2}{2}- \vp\cdot\vlambda\right]
d\xi \ \ \ ; \ \ \forall \xi\in\Lambda_\mu \ . $$
\end{definition}
The existence and uniqueness of weak solution is described in
\rth{th1} below.
\begin{theorem}\label{th1}
For any $\mu\in\overline{\cal M}$  and $\vlambda\in\R^n$,
  there exists a unique
weak solution $\nu\in \Lambda_\mu$ of $\oF( \mu,\vlambda)$.
Moreover, \be\label{defF} \oF(\mu,\vlambda)=-\int_{ \ttt\times\R^n}
\left[\frac{|\vp|^2}{2}-\vp\cdot\vlambda\right] d\nu \ , \ee and
 \be\label{lemma2.1eq}
\oF^*(\mu,\vj)
=\frac{1}{2}\inf_{\nu\in\Lambda_\mu^{\vj}}\int_{\ttt\times\R^n}|\vp|^2d\nu
\ , \ee where the RHS of (\ref{lemma2.1eq}) is attained for the weak
solution of $\oF(\mu,\vlambda)$, provided\footnote{The existence of
a subgradient of $\oF$ with respect to $\vlambda$, among other
results, is stated and proved in Lemma~\ref{fund},
section~\ref{aux}.} $\vj\in\partial_{\vlambda}\oF(\mu,\vlambda)$.
\par
If $\mu(\dx)=\rho(x)\dx$ where $\rho\in C^1(\ttt)$ and $\rho>0$ on
$\ttt$, then a weak solution $\nu$  of $\oF(\mu,\vlambda)$ takes the
form $\nu(\dx d\vp)
=\delta^{\vp}_{\vlambda+\nabla\phi(x)}\otimes\rho(x)dx$ where $\phi$
is the classical solution of the elliptic equation.
\be\label{classweakeq} \nabla\left[
\rho(\nabla_x\phi+\vlambda)\right]=0 \ \ \ \text{on} \ \mathbb{T}^n
\ . \ee
\end{theorem}
\begin{remark}
Equation (\ref{classweakeq}) is strongly elliptic equation if
$\rho>0$, so it has a unique (up to a constant), classical solution.
See, e.g. [GT].
\end{remark}
\begin{remark}
As a by-product we obtain the relation $$\oF^*(\mu,\vj)=
\inf_{\phi\in C^1(\ttt)} {\cal E}(\mu,\vj,\phi) \ , $$ see
\rlemma{fund}-(3).
\end{remark}

\noindent {\bf Example}~1: If $\mu(\dx)=\alpha \rho(x)\dx +
(1-\alpha) \delta_{x_0}$ then the weak solution associated with
$\vlambda\in\R^n$ is
$$ \nu= \delta^{\vp}_{\vlambda+\nabla\phi}\otimes \alpha\rho(x)dx +
(1-\alpha)\delta^{\vp}_0\otimes\delta_{x_0} \ , $$ where $\phi$ is
the classical solution of (\ref{classweakeq}).
\par
We may observe that $\oF$ is concave on $\oM$ for fixed
$\vlambda\in\R^n$. In particular, it is {\it upper-semi-continuous}
in the $C^*$ topology of $\oM$: \be\label{weakinq}
\lim_{n\rightarrow\infty} \oF(\mu_n,\vlambda) \leq \oF(\mu,\vlambda)
\ee whenever $\mu_n\rightarrow \mu$ in $C^*(\mathbb{T}^n)$.
\par\noindent {\bf Example}~2: \ If $\mu_n$ is an atomic measure
then $\oF(\mu_n, \vlambda)=0$ for any $\vlambda\in\R^n$. In
particular, if $\mu_n= n^{-1}\sum_{j=1}^n \delta_{x^{(n)}_j}$ is a
sequence of empirical measures approximating $\mu\in\oM$ then the
L.H.S  of (\ref{weakinq}) is identically zero.
\par
 In \rth{th3} and
Corollary~\ref{cor1} we demonstrate that, in an appropriate sense,
any weak solution is a limit of classical ones.
\begin{theorem}\label{th3}
If $\lim\mu_j=\mu$ in the $C^*(\ttt)$  and
 \be\label{spapp}
\lim_{j\rightarrow\infty} F(\mu_j, \vlambda)= F(\mu, \vlambda)\ \ \
 \ee
 holds, then there exists a subsequence of weak solutions $\nu_j$ of
$\oF(\mu_j,\vlambda)$ along which
 \be\label{limofnu} \lim_{j\rightarrow\infty}\nu_j =\nu \ee
 holds in $C^*(\ttt\times\R^n)$, where
 $\nu$ is a weak solution of
$\oF(\mu,\vlambda)$. \par Moreover, there exists a sequence of
smooth measures $\mu_j=\rho_j \dx$ so that $\rho_j\in
C^\infty(\ttt)$ and $\rho_j>0$ on $\ttt$, for which (\ref{spapp})
holds for any $\vlambda\in\R^n$.
\end{theorem}
\begin{corollary}\label{cor1}
The weak solution of $\oF(\mu,\vlambda)$ is the weak limit
$$ \lim_{j\rightarrow\infty} \delta_{\vp-\vlambda-\nabla\phi_j}d\vp\otimes\rho_j dx  = \nu$$
  where $\phi_j$ are the solutions of
  $$ \nabla\cdot\left( \rho_j(\nabla\phi_j+\vlambda)\right)=0$$
  and $\mu_j=\rho_j\dx\rightarrow \mu$ as guaranteed by \rth{th3}.
\end{corollary}

\begin{definition}\label{G}
Given a continuous function $\Xi\in C(\ttt)$, \be\label{oG*}
\oE_\Xi^*( \vj ):= \sup _{\mu\in \oM} \int \Xi d\mu- \oF^*(\mu, \vj)
\ . \ee Likewise \be\label{oG} \oE_\Xi( \vlambda ):= \sup _{\mu\in
\oM} \int \Xi d\mu +  \oF(\mu, \vlambda)  \ . \ee
\end{definition}
\begin{lemma}\label{dualkj}
$\oE_\Xi^*$ is the negative of the convex dual of $\oE_\Xi$ with
respect to $\R^n$. Than is: \ \ \ $ \oE_\Xi^*(\vj
)=-\sup_{\vlambda\in\R^n} \left\{ \vlambda\cdot\vj-
\oE_\Xi(\vlambda)\right\}$.
\end{lemma}
\begin{proposition}\label{cordelta2}  \  $\oE_\Xi(\vlambda) \geq
\max_{{\ttt}} \Xi$ and $\oE_\Xi^*(\vj) \geq \max_{{\ttt}} \Xi$ hold
for any $\vlambda, \vj \in\R^n$. If $|\vlambda|$ (res. $|\vj|$) is
large enough, then the inequality is strong.
\end{proposition}
\vskip .2in\noindent {\bf Open problem:} {\it Is $\oE_\Xi(\vlambda)
> \max_{{\ttt}} \Xi$ for any $\vlambda\not=0$ ?  }
\begin{lemma}\label{muexists}
For any $\Xi\in C(\ttt)$  $\vlambda\in\R^n$ there exists
$\mu_0\in\oM$ verifying the maximum in (\ref{oG}).
 There exists $\vj\in\partial_{\vlambda}\oF(\mu_0,\Xi)\subset
\partial_{\vlambda} \oE_\Xi(\cdot)$ \ for which $\mu_0$  verifies the
maximum in (\ref{oG*}).
\end{lemma}
\vskip .2in We end this section by stating the connection between
maximizers of $\oE_\Xi$ and $\oE_\Xi^*$ and the minimal invariant
measures of the weak-KAM theory:

\begin{theorem}\label{connect}
If \ $\Xi$ is smooth enough (say, $\Xi\in C^2(\ttt)$) and $\nu$ is a
Mather measure (\ref{minmeasL}) of the Lagrangian $
L=|\vp|^2/2-\Xi(x)-\vp\cdot\vlambda$ on $\mathbb{T}^n\times \R^n$
and $\mu=\pi_\#\nu$ then $\nu$ is weak solution of $\oF(\mu,
\vlambda)$, and is a maximizer of $\oE_\Xi(\vlambda)$ in (\ref{oG}).
\par
Moreover, $\nu$ verifies (\ref{minmeas}) where $\vj$ is the rotation
number   $\vec{\alpha}(\nu)$ given by (\ref{rotnum}).
\end{theorem}

\subsection{Extension to time dependent measures}\label{TD}
We now extend the definition of $\oF$ and $\oF^*$ to the set of
$\oM-$valued  orbits on the interval $[0,T]$.

 Define,
for\footnote{See point \ref{orbits} in section~\ref{list} . }
$\hat{\mu}\in \oM_T$  \be\label{FbarT} F(\hat{\mu},\vlambda, T)
:=\frac{1}{T}\inf_{\phi\in
\CT}\int_0^T\int_{{\ttt}}\left[\frac{\partial\phi}{\partial t}+
\frac{1}{2}|\nabla_x\phi+\vlambda|^2 \right]d\mu_{(t)} dt \ . \ee
Let also ${\cal E}$ a function on $\oM_T\times\R^n\times \CT$
defined as:
 \be\label{calET} {\cal E}(\hat{\mu},\vj,\phi,T):= \left\{\frac{1}{2}
\left| \vj-\frac{1}{T}\int_0^T\int_{{\ttt}}\nabla\phi
d\mu_{(t)}dt\right|^2- \frac{1}{T}\int_0^T\int_{{\ttt}}\left(
\phi_t+\frac{1}{2}|\nabla\phi|^2\right)d\mu_{(t)} dt \right\} \ .
\ee

 The analog
of Definition~\ref{deflift} is
\begin{definition}:\label{def4.5}
\begin{enumerate}
\item
The set $\widehat{\Lambda}_T \subset \M_T$ consists of all orbits of
probability measures $\hat{\nu}:[0,T]\rightarrow{\cal
M}(\mathbb{T}^n\times\R^n)$, $\left. \hat{\nu}\right|_{(t)}:=
\nu_{(t)}\in {\cal M}(\mathbb{T}^n\times\R^n)$, for which
$$\int_0^T\int_{ \ttt\times\R^n }\left(\theta_t+ \vp\cdot\nabla_x\theta(x,t)\right) d\nu_{(t)}dt = 0 \ \ \
\forall\theta\in C_0^1(\ttt\times (0,T)) \text{and}  \int_0^T\int_{
\ttt\times\R} |\vp|^2d\nu_{(t)}dt<\infty \ .
$$
\item
Given  $\hat{\mu}\in \oM_T$, the set $\widehat{\Lambda}_{T,
\hat{\mu}}\subset \Lambda$ of liftings of $\hat{\mu}$ is composed of
all $\hat{\nu}\in\widehat{\Lambda}_T$ for which
$\pi_{\#}\hat{\nu}=\hat{\mu}$, that is,
$\pi_{\#}\nu_{(t)}=\mu_{(t)}$ for any $t\in[0,T]$.
\item For each $\vj\in\R^n$, the set $\widehat{\Lambda}^{\vj}_{T, \hat{\mu}}$ is
defined as all  $\hat{\nu}\in\widehat{\Lambda}_{T,\hat{\mu}}$ which
satisfies
$$ T^{-1}\int_0^T\int_{ \ttt\times\R^n} \vp d\nu_{(t)} dt = \vj \ . $$
\end{enumerate}
\end{definition}
\begin{definition}\label{defweakT}
For given $\vlambda\in\R^n$,
$\hat{\nu}\in\widehat{\Lambda}_{T,\hat{\mu}}$ is
 a weak solution of
$\oF(\hat{\mu},\vlambda,T)$  provided
$$ \int_0^T\int_{ \ttt\times\R^n} \left[\frac{|\vp|^2}{2}-\vp\cdot\vlambda\right] d\nu_{(t)}dt \leq \int_0^T\int_{ \ttt\times\R^n}
 \left[\frac{|\vp|^2}{2}- \vp\cdot\vlambda\right]
d\xi_{(t)} \ \ \ ; \ \ \forall
\hat{\xi}\in\widehat{\Lambda}_{\hat{\mu}} \ .
$$
\end{definition}\vskip .2in \noindent
The $T-$orbit analogue of \rth{th1} is
\begin{proposition}\label{thT}
For any $\hat{\mu}\in\overline{\cal M}_T$  and $\vlambda\in\R^n$,
  there exists a unique
weak solution $\hat{\nu}^{(0)}\in \widehat{\Lambda}_{T,\hat{\mu}}$
of $\oF( \hat{\mu},\vlambda,T)$. Moreover, \be\label{defFT}
\oF(\hat{\mu},\vlambda,T)=-\frac{1}{T}\int_0^T\int_{\ttt\times\R^n}
\left[\frac{|\vp|^2}{2}-\vp\cdot\vlambda\right] d\nu^{(0)}_{(t)} dt
=-
\inf_{\hat{\nu}\in\widehat{\Lambda}_{T,\hat{\mu}}}\frac{1}{T}\int_0^T\int_{\ttt\times\R^n}
\left[\frac{|\vp|^2}{2}-\vp\cdot\vlambda\right] d\nu_{(t)} dt
 \ . \ee The Legendre transform of $\oF(\hat{\mu}, \vlambda,T)$
with respect to $\vlambda$  is \be\label{lemma2.1eqT}
\oF^*(\hat{\mu},\vj,T)=\sup_{\phi\in \CT}{\cal
E}(\hat{\mu},\vj,\phi)
=\frac{1}{T}\inf_{\hat{\nu}\in\widehat{\Lambda}_{T,\hat{\mu}}^{\vj}}\int_0^T\int_{\ttt\times\R^n}\left[\frac{|\vp|^2}{2}-\vp\cdot\vlambda\right]d\nu_{(t)}dt
\ , \ee where the RHS of (\ref{lemma2.1eqT}) is attained for the
weak solution of $\oF(\hat{\mu},\vlambda,T)$, provided
$\vj\in\partial_{\vlambda}\oF(\hat{\mu},\vlambda,T)$.
\par
If $\mu_{(t)}(\dx)=\rho(x,t)\dx$ where $\rho\in \CT$ and $\rho>0$ on
$\ttt\times [0,T]$ then a weak solution $\hat{\nu}$ of
$\oF(\hat{\mu},\vlambda,T)$ takes the form $\nu_{(t)} =
\delta_{\vp-\vlambda-\nabla\phi(x,t)} d\vp\otimes \rho(x,t)dx$ where
$\phi$ is the classical solution of the elliptic equation.
\be\label{classweakeqT} \nabla_x\left[
\rho(\nabla_x\phi+\vlambda)\right]=-\rho_t  \ . \ee
\end{proposition}

\begin{definition}\label{DefFT}Given $\mu_1, \mu_2 \in\oM$, $\vlambda\in\R^n$ and $\Xi\in C(\ttt)$,
define\footnote{Recall section~\ref{list}-(\ref{orbits}) .}
\be\label{FbarTmin} H_{\Xi,T}(\mu_1, \mu_2,\vlambda):=
\sup_{\hat{\mu}\in \oM_T(\mu_1,\mu_2)} \left[ \frac{1}{T}\int_0^T
\int_{{\ttt}}\Xi d\mu_{(t)}dt+ \oF(\hat{\mu}, \vlambda, T)\right] \
. \ee Likewise, for any $\vj\in\R^n$: \be\label{FbarTmin*}
H^*_{\Xi,T}(\mu_1, \mu_2,\vj):= \sup_{\hat{\mu}\in
\oM_T(\mu_1,\mu_2)}\left[\frac{1}{T}\int_0^T \int_{{\ttt}}\Xi
d\mu_{(t)}dt- \oF^*(\hat{\mu}, \vj, T) \right] \ . \ee If
$\mu_1=\mu_2\equiv \mu\in \oM$ then $H_{\Xi,T}(\mu,\vlambda
):=H_{\Xi,T}(\mu,\mu,\vlambda)$  and $H^*_{\Xi,T}(\mu,\vj
):=H^*_{\Xi,T}(\mu,\mu, \vj)$.  \par\noindent
 If $\Xi\equiv 0$, set
$\oF_T(\mu,\vlambda):=H_{0,T}(\mu,\vlambda)$ and $F^*_T(\mu,\vj):=
-H^*_{0,T}(\mu,\vj)$.
\end{definition}

\begin{proposition}\label{maxoverline}
For any $\Xi\in C(\ttt)$, $\vlambda\in\R^n$  and $\mu_1,\mu_2\in\oM$
there exists an orbit $\hat{\mu}\in\oM(\mu_1,\mu_2)$ realizing
(\ref{FbarTmin}). Likewise, for any $\vj\in\R^n$  there exists an
orbit $\hat{\mu}\in\oM(\mu_1,\mu_2)$ realizing (\ref{FbarTmin*}).
\end{proposition}
\begin{proposition}\label{prop4.3}
For any  $T>0$, $\Xi\in C(\ttt)$,  $\vlambda\in\R^n$ and
$\mu\in\oM$, \be\label{obvineq} H_{\Xi,T}(\mu,\vlambda)\geq
\int_{{\ttt}}\Xi d\mu+ \oF_T(\mu,\vlambda) \geq \int_{{\ttt}}\Xi
d\mu+ \oF(\mu,\vlambda) \ , \ee but \be\label{oF=oG}\sup_{\mu\in
\oM} H_{\Xi,T}(\mu,\vlambda)=
\oE_\Xi(\vlambda):=\sup_{\mu\in\oM}\int_{{\ttt}}\Xi d\mu+
\oF(\mu,\vlambda)= \sup_{\mu\in\oM}\int_{{\ttt}}\Xi d\mu+
\oF_T(\mu,\vlambda) \ \ \forall \ T>0 \ , \ee and the maximizer of
$\oE_\Xi(\vlambda)$ (\ref{oG})  is the same as the maximizer of
$H_{\Xi,T}(\mu,\vlambda)$ for any $T>0$.
 \end{proposition}

 \subsection{Optimal transportation}
\label{MK}
 Recall the
definition of the action associated with the Lagrangian
(\ref{mechkag}):  \\ $A^\Xi_{\vlambda}:  \ttt\times \ttt \times
\R\rightarrow \R$: \be\label{Adef} A^\Xi_{\vlambda}(\vy, \vx,
T)=\inf_{\vx(\cdot)}\left\{\frac{1}{T}\int_0^T \left(
\frac{|\vdotx-\vlambda |^2}{2} - \Xi(x(s))\right)ds \ \ \ ; \ \ \
\vx(0)=\vy \ , \vx(t)=\vx \ ,\right\} \ee
\begin{definition} For
$\vlambda\in\R^n$, $\mu_1,\mu_2\in\oM$
 define\footnote{
Recall section~\ref{list}-(\ref{listm^2}).} the Monge-Kantorovich
distance with respect to the action $A^\Xi$:  \be\label{dvlambda}
D^T_{\vlambda}(\mu_1,\mu_2, \Xi):= \min_{\sigma\in \oM^{(2)}(\mu_1,
\mu_2)} \int_{ \ttt}A^\Xi_{\vlambda}(x,y,T)d\sigma(x ,y) \ .  \ee If
$\mu_1=\mu_2\equiv \mu\in\oM$ we denote \be\label{dvlambda1}
D^T_{\vlambda}(\mu, \Xi):= \min_{\sigma\in \oM^{(2)}(\mu)} \int_{
\ttt}A^\Xi_{\vlambda}(x,y,T)d\sigma(x ,y) \ .  \ee
\end{definition}
\begin{example}\label{MKXi=0} If $\Xi\equiv 0$ then $A^0_{\vlambda}(\vy,\vx,T)=
\|\vx-\vy-T\vlambda\|^2/ (2T^2)$, where $\| \cdot \|$ is the
Euclidian metric on $\ttt$. In particular \be\label{DTk0}
D^T_{\vlambda} (\mu):= D^T_{\vlambda}(\mu, 0)= \min_{\sigma\in
\oM^{(2)}(\mu)}\frac{1}{2T^2} \int_{ \ttt}\| x-y-T\vlambda\|^2
d\sigma(x, y) \ . \ee
\end{example}

\begin{proposition}\label{FtoD1}
For any $\mu_1, \mu_2\in\oM$, $\Xi\in C(\ttt)$, $\vlambda\in\R$ and
$T>0$
$$ H_{\Xi,T}(\mu_1, \mu_2,\vlambda)= \frac{|\vlambda|^2}{2}-
D_{\vlambda}^T(\mu_1, \mu_2, \Xi)\ \ .  $$ In particular
\be\label{BBthC} H_{\Xi,T}(\mu,\vlambda)= \frac{|\vlambda|^2}{2}-
D_{\vlambda}^T(\mu, \Xi)\ \ .  \ee
\end{proposition}
\begin{lemma}
\label{masstrans} $D^T_{\vlambda}: \oM\times C(\ttt)\rightarrow \R$
is continuous in both the weak $C^*$ topology of $\oM$ and in the
$\sup$ topology of $C(\ttt)$.
\end{lemma}
\begin{proposition}\label{FtoD2} For any $\mu\in\oM$, $\Xi\in C(\ttt)$, $\vlambda\in\R$
 \be\label{limT}
\lim_{T\rightarrow 0} H_{\Xi,T}(\mu,\vlambda)= \int_{{\ttt}}\Xi d\mu
+ \oF(\mu,\vlambda) \ . \ee
\end{proposition}

As a corollary to (\ref{limT}) and Lemma~\ref{masstrans}, evaluated
for $\Xi=0$,  we obtain
\begin{theorem}\label{limtt}
The functional $H_{\Xi,T}$ (Definition~\ref{DefFT}) is continuous on
$\oM$ with respect to the $C^*$ topology. In addition
$$ \lim_{T\rightarrow 0} H_{\Xi,T}(\mu,\vlambda)=H_\Xi(\mu,\vlambda)$$
for any $\mu\in\oM$, $\vlambda\in\R^n$.
\end{theorem}
By definition~\ref{DefFT} and Theorem~\ref{limtt}
\begin{corollary}\label{corFT0}
$$ \lim_{T\rightarrow 0} F_T(\mu,\vlambda)=F(\mu,\vlambda)$$
for any $\mu\in\oM$, $\vlambda\in\R^n$.
\end{corollary}

\subsection{A combinatorial search algorithm}

Next, Birkhoff's Theorem implies
\begin{lemma}\label{MKDis}Given $x_1, \ldots x_j\in \ttt$. For $\mu_j:=j^{-1}\sum_1^j
\delta_{x_i}$,
$$ D^T_{\vlambda}(\mu_j, \Xi) = \min_{\sigma\in \Pi_j} \sum_{i=1}^j
A_{\vlambda}^\Xi(x_i, x_{\sigma(i)}, T)$$ where $\Pi_j$ is the set
of all permutations of $\{1, \ldots j\}$. In particular
$$ D^T_{\vlambda}(\mu_j) =\frac{1}{2T^2} \min_{\sigma\in \Pi_j} \sum_{i=1}^j
\| x_i-x_{\sigma(i)}-T\vlambda\|^2 \ . $$

\end{lemma}
By Proposition~\ref{FtoD1}, Theorem~\ref{limtt},
Corollary~\ref{corFT0} and Lemma~\ref{MKDis} we obtain the following
algorithm for evaluation of $H_\Xi(\mu, \vlambda)$ and $F(\mu,
\vlambda)$:
\begin{theorem}\label{disFH}
 Let $\mu_j:= j^{-1}\sum_{i=1}^j
\delta_{x^{(j)}_i}$ is a sequence of measures  converging $C^*$ to
$\mu\in \oM$. Then \be\label{say3} \lim_{T\rightarrow 0}
\lim_{j\rightarrow\infty}
D^T_{\vlambda}(\mu_j,\Xi)=\frac{|\vlambda|^2}{2}- H_\Xi(\mu,
\vlambda) \ . \ee In particular
$$ \lim_{T\rightarrow 0}
\lim_{j\rightarrow\infty} \frac{1}{2T^2} \min_{\sigma\in \Pi_j}
\sum_{i=1}^j \|
x^j_i-x^j_{\sigma(i)}-T\vlambda\|^2=\frac{|\vlambda|^2}{2}-\oF(\mu,
\vlambda) \ .
$$

\end{theorem}
We may use now Theorem~\ref{disFH} to evaluate the effective
Hamiltonian $\oE_\Xi(\vlambda)$. In fact, we do not need to take the
limit $T\rightarrow 0$, as shown below:

\begin{definition}
Given $j\in \mathbb{N}$, let
$$D^T_{\vlambda}(j, \Xi):= \min_{x_1, \ldots
x_j\in\ttt}\min_{\sigma\in \Pi_j} \sum_{i=1}^j A_{\vlambda}^\Xi(x_i,
x_{\sigma(i)}, T) \  $$ where $\Pi_j$ as defined in
Lemma~\ref{MKDis}. Let also \be\label{barDpdef}
\overline{D}^T_{\vlambda}(j,\Xi):= \min_{x_1, \ldots
x_j\in\ttt}\min_{\sigma\in \Pi_j}\sum_{l=1}^j \left[(2T^2)^{-1}
\|x^{(j)}_l-x^{(j)}_{\sigma(l)} -
 T\vlambda\|^2+ j^{-1} \Xi(x_i)\right] \ \ .
 \ee
 In particular, for $\Xi=0$,
$$D^T_{\vlambda}(j, 0)=\overline{D}^T_{\vlambda}(j,0):=
D^T_{\vlambda}(j) \ . $$

\end{definition}
\begin{theorem}\label{DisH}
For any $\Xi\in C(\ttt)$, $\vlambda\in\R^n$ and $T>0$,
\be\label{athdish}\oE_\Xi(\vlambda)= \frac{|\vlambda|^2}{2}-
\lim_{j\rightarrow\infty} D^T_{\vlambda}(j, \Xi) \ .  \ee
\end{theorem}

\section{Proof of the main results}\label{proofs}
\subsection{Duality}
 The key duality argument for
minimizing convex functionals under affine constraints  is
summarized in the following proposition. This is a slight
generalization of Proposition~4.1 in [W1]. The proof  is sketched in
the appendix of [W1]).
\begin{proposition}\label{propdual}
Let ${\bf C}$ a real Banach space and  ${\bf C}^*$ its dual. Let
${\bf Z}$ a   subspaces of \ ${\bf C}$. Let $h\in {\bf C}^*$. Let
${\bf Z}^*\subset {\bf C}^*$ given by the condition $z^*\in {\bf
Z}^*$ iff $<z^*-h,z>=0$ for any $z\in {\bf Z}$.   Let ${\cal F}:{\bf
C}^*\rightarrow \R\cup\{\infty\}$ a convex function and
\be\label{Edefth1} E:= \inf_{c^*\in {\bf Z}^*} {\cal F}(c^*)  \ .
\ee Then \be\label{Edef2th1} \sup_{z\in {\bf Z}}\inf_{c^*\in {\bf
C}^*} \left[ {\cal F}(c^*)- <c^*-h,z> \right]\leq E \ , \ee
 and if   $\overline{A}_0:= \{ c^*\in {\bf C}^* \ ; \ {\cal F}(c^*) \leq
E\}$ is compact (in the $*-$ topology of ${\bf C}^*$), then there is
an equality in (\ref{Edef2th1}).

 In particular, $E<\infty$ if  and only if
${\bf Z}^*\not=\emptyset$. In this case there exists $z^*\in {\bf
Z}^*$ for which $E={\cal F}(z^*)$.
\end{proposition}
\begin{remark}
The case $E<\infty$ does not implies, in general, the existence of
$z\in{\bf Z}$ realizing (\ref{Edef2th1}) .
\end{remark}
\subsection{An Auxiliary  result}\label{aux}
\begin{lemma}\label{fund} \ .
\begin{enumerate}
\item
$\oF$ is  convex on $\R^n$ as function of $\vlambda$ and   concave
on $\oM$ as function of $\mu$.
\item $\oF^*$ is convex on both $\R^n$ (as a function of $\vj$) and
$\oM$.
\item \ \ \ \ \ $\oF^*(\mu,\vj)= \inf_{\phi\in C^1(\ttt)} {\cal
E}(\mu,\vj,\phi)$ \ .
\item
The sub-gradients $\partial_{\vlambda} \oF(\mu, \cdot)$ and
$\partial_{\vj} \oF^*(\mu, \cdot)$ are non-empty convex cones in
$\R^n$  for any $\mu\in\oM$ and $\vlambda\in\R^n$ (res.
$\vj\in\R^n$) and satisfies $\vlambda\in\partial_{\vj} \oF^*(\mu,
\cdot)$ iff $\vj\in\partial_{\vlambda} \oF(\mu, \cdot)$.
\item\label{usc}
 $\oF$ is upper-semi-continuous in
the $C^*$ topology of  $\oM$ for any $\vlambda\in\R^n$, and $\oF^*$
is {\it lower}-semi-continuous for the same topology for nay
$\vj\in\R^n$.
\end{enumerate}
\end{lemma}
In particular, from point (5) of this Lemma:
\begin{corollary}\label{usefordelta}
 If $x_0\in \ttt$ and  $\mu_n\rightarrow \delta_{x_0}$ then
$\lim_{n\rightarrow \infty} \oF(\mu_n, \vlambda)= \oF(\delta_{x_0},
\vlambda)=0$.
\end{corollary}
\noindent {\bf Proof of Lemma~\ref{fund}}:  \ \ Concavity of $\oF$
on $\oM$ is a result of its definition as a n infimum of functionals
over this convex set. The strict convexity on $\R^n$ follows from
its quadratic dependence:
$$ \frac{|\vlambda_1+\nabla\phi_1|^2}{2} +\frac{|\vlambda_2+\nabla\phi_2|^2}{2}\geq
\left|\frac{\vlambda_1+\vlambda_2}{2}+\frac{\nabla\phi_1+\nabla\phi_2}{2}\right|^2
 \ , $$
holds for any $\vlambda_1,\vlambda_2\in\R^n$ and any
$\phi_1,\phi_2\in C^1(\ttt)$.  If $\phi_1$ approximates the
maximizer of $\oF(\mu,\vlambda_1)$ (res. $\phi_2$ approximates the
maximizer of $\oF(\mu,\vlambda_2)$), then integrating the above
inequality with respect to $\mu$ yields
$$ \oF(\mu, \vlambda_1) + \oF(\mu, \vlambda_2)\geq
\int_{{\ttt}}\left|\frac{\vlambda_1+\vlambda_2}{2}+\frac{\nabla\phi_1+\nabla\phi_2}{2}\right|^2d
\mu \geq  2\oF\left(\mu, \frac{\vlambda_1+\vlambda_2}{2}\right) \ .
$$
The same arguments apply to $\oF^*$. In addition, from (\ref{Fbar2},
\ref{F*bar2}) $$\oF^*(\mu,\vj)= \sup_{\vlambda\in \R^n} \left[
\vlambda\cdot\vj - \oF(\mu,\vlambda)\right] = \sup_{\phi\in
C^1(\ttt)} \left\{ \sup_{\vlambda\in \R^n} \left[ \vlambda\cdot\vj -
\frac{1}{2}\int_{{\ttt}}|\vlambda+\nabla\phi|^2 d\mu \right]\right\}
\ , $$ and from $$ {\cal E}(\mu,\vj,\phi)=\sup_{\vlambda\in \R^n}
\left[ \vlambda\cdot\vj -
\frac{1}{2}\int_{{\ttt}}|\vlambda+\nabla\phi|^2 d\mu \right] \  $$
we obtain (3). (4-6) follow from (1,2).
\subsection{Proof of \rth{th1}}
 We now apply Proposition~\ref{propdual} as follows:

Let ${\bf C}$ the space of all continuous functions on
$\ttt\times\R^n$, equipped with the norm
$$\| q\|:=  \sup_{(x,\vp)\in \ttt\times
\R^n}\left\{\frac{|q(x,\vp)|}{1+|\vp|} \right\}  \ . $$ Define
\be\label{Zdef} {\bf Z}:= \left\{ \vp\cdot\nabla_x\phi  \ \ , \ \
\phi\in C^1(\ttt) \  \right\} \ee The dual space ${\bf C}^*$ {\it
contains} the set ${\cal M}(\ttt\times\R^n)$ of finite Borel
measures on $\ttt\times\R^n$ which admit a finite first moment. If
$\nu\in{\bf C}^*$ is such a measure then the duality relation is
given by
$$ \langle\nu, q\rangle= \int_{
\ttt\times\R^n}qd\nu  \  \forall q\in {\bf C} \ . $$ Given a
probability measure $\mu\in \oM$, define
$$ {\cal F}_\mu(\nu) =\left \{ \begin{array}{c}
 \int_{\ttt\times\R^n} \left( |\vp|^2/2 -\vlambda\cdot\vp\right)d\nu  \ \ \text{if}
  \ \ \nu\in {\cal M}(\ttt\times\R^n)\cap {\bf C}^* \
    \text{and satisfies} \ \nu(\dx, \R^n) = \mu(\dx)  \\
  \infty \ \ \text{otherwise} \\
\end{array}\right. $$
(recall Definition~\ref{deflift}).
 Evidently, ${\cal F}_\mu$ is a convex function on ${\bf C}^*$. Note
 also that the set $\overline{A}_0:= \{ c^* \ ; {\cal F}_\mu(c^*)<E\}\subset {\bf C}^*$ is compact
 for any $E<\infty$ by Prokhorov Theorem.
\begin{lemma}\label{lem3.1}
 $\nu\in {\bf Z}^*$ and ${\cal
F}_\mu(\nu) < \infty$ if and only if $\nu\in \Lambda_\mu$.
\end{lemma}
Substitute this ${\cal F}_\mu$ for ${\cal F}$ in (\ref{Edefth1})
where $h\equiv 0$  it follows that \be\label{E=} E =
\inf_{\nu\in\Lambda_\mu}\int_{\ttt\times\R^n}
\left(|\vp|^2/2-\vlambda\cdot\vp\right) \nu(\dx)
 \ .
\ee On the other hand \begin{multline}\label{111} \inf_{c^*\in {\bf
C}^*}{\cal F}(c^*)- <c^*,z> = \inf_{\nu\in {\bf
C}^*}\int_{\ttt\times\R^n} \left( \frac{1}{2}|\vp|^2- \vp\cdot
(\nabla_x\phi +\vlambda) \right)d\nu \\
=\inf_{\nu\in\Lambda_\mu}\frac{1}{2}\int_{\ttt\times\R^n}
|\vp-\vlambda-\nabla_x\phi|^2d\nu-\frac{1}{2}\int_{{\ttt}}|\nabla_x\phi+\vlambda|^2d\mu
 \ .
\end{multline}
We choose $\nu= \delta^p_{(\vlambda+\nabla_x\phi)}\otimes\mu$, so
the first term in (\ref{111}) is zero. Then \be\label{legendreoF}
\sup_{z\in{\bf Z}}\inf_{c^*\in {\bf C}^*}{\cal F}(c^*)-
<c^*,z>=\sup_{  \phi\in C^1(\ttt)}\left[
-\frac{1}{2}\int_{{\ttt}}|\nabla_x\phi+\vlambda|^2d\mu\right]:=
-\oF(\mu,\vlambda) \ .
 \ee
where $\oF$ as defined in (\ref{Fbar2}). The last part of
Proposition~\ref{propdual} implies the existence of a weak solution
$\nu\in\Lambda_\mu$ of $\oF(\mu,\vlambda)$. \par To show the
uniqueness of the weak solution, note that any Borel measure $\nu$
on $\ttt\times\R^n$ whose $\ttt$ marginal is $\mu$ can be written as
$\nu(\dx d\vp) = \mu(\dx)Q_{x}(d\vp)$ where $Q_x$ is a Borel
probability measure on $R^n$ defined for $\mu-$a.e. $x\in \ttt$. If
$\nu$ satisfies Definition~\ref{defweak}, then
$Q_x=\delta^p_{\vec{v}}$ where $\vec{v}(x):= \int_{\R^n} \vp
Q_x(d\vp)$ is a Borel vector filed, defined $\mu-$a.e. If there are
$\nu_1\not=\nu_2$ which realize the minimum in
Definition~\ref{defweak} and $\vec{v}_1, \vec{v}_2$ the
corresponding vector fields, then
$$ \oF(\mu,\vlambda)=\frac{1}{2}\int_{{\ttt}}\left|\vec{v}_1\right|^2d\mu=
\frac{1}{2}\int_{{\ttt}}\left|\vec{v}_2\right|^2d\mu$$ implies
$$\frac{1}{2} \int_{{\ttt}}\left|\frac{\vec{v}_1+\vec{v}_2}{2}\right|^2d\mu  <\oF(\mu,\vlambda) \ , $$
unless $\vec{v}_1=\vec{v}_2$ $\mu-$ a.e., which contradicts  the
minimality of $\nu_1$ and $\nu_2$. \par From  (\ref{legendreoF}) it
follows that the {\it Legendre transform} of the function $\oF(\mu,
\cdot)$ is
 \begin{multline}\label{cale}\oF^*(\mu,\vj)=
\sup_{\vlambda\in\R^n} \sup_{\phi\in C^1(\ttt)}\left[
\vlambda\cdot\vj-\frac{1}{2}\int_{{\ttt}}|\nabla_x\phi+\vlambda|^2d\mu\right]
\\ =  \sup_{\phi\in C^1(\ttt)}\sup_{\vlambda\in\R^n} \left[
\vlambda\cdot\vj-\frac{1}{2}\int_{{\ttt}}|\nabla_x\phi+\vlambda|^2d\mu\right]
=\sup_{\phi\in C^1(\ttt)} {\cal E}(\mu,\vj,\phi)
\end{multline} where ${\cal E}$ as defined in (\ref{calE}).
\par
 To prove the last part, note that
$\nu_0:=\rho(x)dx\otimes \delta_{(\vp-\nabla\psi-\vlambda)} d\vp\in
\Lambda_\mu$ whenever $\psi$ is the solution of (\ref{classweakeq}).
Indeed, by (\ref{classweakeq}) and integration by parts
\be\label{www}\int_{\ttt\times\R^n}\nabla\phi\cdot\vp d\nu_0=
\int_{{\ttt}} \nabla\phi\cdot(\vlambda+\nabla\psi)\rho(x) \dx =
\int_{{\ttt}}
\phi\nabla\cdot\left[(\vlambda+\nabla\psi)\rho(x)\right] \dx  =0 \
\ee for any $\phi\in C^1(\ttt)$. Hence, (\ref{E=}) implies
\be\label{Eleq} E \leq \int_{\ttt\times\R^n}
\left(|\vp|^2/2-\vlambda\cdot\vp\right) d\nu_0 =
\frac{1}{2}\int_{{\ttt}}\left| \nabla\psi + \vlambda\right|^2
\rho(x)\dx - \vlambda\cdot \int_{{\ttt}} \left(
\nabla\psi+\vlambda\right) \rho(x)\dx \ . \ee However, $\psi$
realizes the infimum in (\ref{Fbar2}), so
(\ref{E=},\ref{legendreoF}) and Proposition~\ref{propdual} imply
$$ E=-\frac{1}{2}\int_{{\ttt}}\left| \nabla\psi + \vlambda\right|^2
\rho(x)\dx$$ which, together with (\ref{Eleq}), imply
$$\int_{{\ttt}}\left| \nabla\psi + \vlambda\right|^2
\rho(x)\dx - \vlambda\cdot \int_{{\ttt}} \left(
\nabla\psi+\vlambda\right) \rho(x)\dx\geq 0 \ . $$ However,
(\ref{www}) with $\phi=\psi$ implies the equality above, hence the
equality in (\ref{Eleq}) as well. In particular, $\nu_0$ minimizes
(\ref{E=}).
  \ \ \ \ \ $\Box$
\subsection{Proof of theorem~\ref{th3}} First, the sequence
$\{\nu_j\}$ is tight in ${\cal M}(\ttt\times\R^n)$ since $\ttt$ is
compact and $\int_{\ttt\times\R^n}|\vp|^2 d\nu_j$ are uniformly
bounded. By Prokhorov Theorem it follows that the weak limit
$\nu_j\rightarrow \nu\in{\cal M}(\ttt\times\R^n)$ exists (for a
subsequence). Also, $\nu\in\Lambda_\mu$ since the condition given in
Definition~\ref{deflift} is preserved under the weak-* convergence.
\par
Next
\begin{multline}\label{chain}-\lim_{j\rightarrow\infty}\oF(\mu_j,\vlambda)=\lim_{j\rightarrow\infty}
\int_{ \ttt\times\R^n}\left( |\vp|^2/2-\vlambda\cdot\vp\right)
d\nu_j\geq\int_{ \ttt\times\R^n}\left(
|\vp|^2/2-\vlambda\cdot\vp\right) d\nu  \\ \geq
\inf_{\xi\in\Lambda_\mu}\int_{ \ttt\times\R^n}\left(
|\vp|^2/2-\vlambda\cdot\vp\right) d\xi =- \oF(\mu,\vlambda) \ ,
\end{multline} By assumption (\ref{spapp}) it follows that the
equality folds in (\ref{chain}). In particular, $\nu$ is the weak
solution of $\oF(\mu,\vlambda)$.
 \par
To prove the second part, let $\eta_\eps\in C^\infty(\ttt)$ a
sequence of {\it positive} mollifiers on $\ttt$ satisfying
$\lim_{\eps\rightarrow 0} \eta_\eps=\delta_{(\cdot)}$, and
$\mu_\eps=\eta_\eps *\mu$. Then $\mu_\eps(\dx)=\rho_\eps(x)\dx$
where $\rho_\eps\in C^\infty(\ttt)$ are strictly positive on $\ttt$
and
$$ \lim_{\eps\rightarrow 0} \mu_\eps=\mu \ . $$
Next, let $\nu$ be a weak solution of $\oF(\mu,\vlambda)$ and
$\nu_\eps=\eta_\eps*\nu$. If $\nu=\mu(dx)\nu_x(d\vp)$ and $q=q(p)$
any $\nu$ measurable function, then $\tilde{q}(x):=\int_{\R^n}
q(\vp)\nu_x(d\vp)$ is $\mu$ measurable and   $$\int_{
\ttt\times\R^n}q(\vp)d\nu= \int_{\ttt} \tilde{q}(x)\mu(dx)$$ while
 \begin{multline}\label{qeq} \int_{ \ttt\times\R^n}q(\vp)d\nu_\eps= \int_{
\ttt\times\R^n}\int_{\ttt}dx\eta_\eps(|y-x|)\mu(dy)\nu_y(d\vp)q(\vp)\\
=\int_{\ttt}\int_{\ttt} \eta_\eps(|x-y|)
dx\mu(dy)\tilde{q}(y)=\int_{\ttt} dx\mu(dy)\tilde{q}(y)=\int_{
\ttt\times\R^n}q(\vp)d\nu
\end{multline} for any $\nu-$measurable function $q$ on $\R^n$. Then
$$ \int_{ \ttt\times\R^n} \nabla\theta\cdot\vp d\nu_\eps(x)=
\int_{ \ttt\times\R^n} \nabla(\eta_\eps *\theta)\cdot\vp d\nu(x)=0
$$ for any $\theta\in C^\infty(\ttt)$ hence
$\nu_\eps\in\Lambda_{\mu_\eps}$. \par Define
$$ \vec{v}_\eps(x)=\rho_\eps^{-1}(x)\int_{\R^n}\vp d\nu_\eps(x,d\vp)  \ \ \ , x\in \ttt \ .  $$
Then $\vec{v}_\eps\in C^\infty(\ttt)$ and $\nabla\cdot(\rho_\eps
\vec{v}_\eps)=0$.   Let $\phi_\eps$ be the unique solution of the
elliptic equation
\be\label{all}\nabla\cdot\left(\rho_\eps[\nabla\phi_\eps+\vlambda]\right)=0
\  \  \ \ . \ee   Define
$$ \widehat{\nu}_\eps(\dx d\vp):=
\rho_\eps(x)dx\otimes \delta_{(\vp-\vlambda-\nabla\phi_\eps)}  d\vp
\ .
$$
 Then (\ref{all}) implies, as in (\ref{www}), that
 $\widehat{\nu}_\eps\in\Lambda_{\mu_\eps}$.
 By (\ref{defF}) and (\ref{qeq}) and the second part of
 \rth{th1}
\begin{multline}\label{Fsup*} -\oF(\mu,\vlambda):= \int_{
\ttt\times\R^n}\left(|\vp|^2/2-\vlambda\cdot\vp\right) d\nu= \int_{
\ttt\times\R^n}\left(|\vp|^2/2-\vlambda\cdot\vp\right)d\nu_\eps
\\ \geq \int_{
\ttt\times\R^n}\left(|\vp|^2/2-\vlambda\cdot\vp\right)d\widehat{\nu}_\eps
= -\oF(\mu_\eps,\vlambda) \ , \end{multline} hence \be\label{sof1}
\lim_{\eps\rightarrow 0}\oF(\mu_\eps, \vlambda) \geq \oF(\mu,
\vlambda) \ \ee
    for any $\vlambda\in\R^n$. But, since $\oF$ is concave in $\mu$
    via Lemma~\ref{fund}  it follows that there is an equality in
    (\ref{sof1}).
  \ \ \ $\Box$
  \subsection{Proof of Lemma~\ref{dualkj}} Since $\oF$ is
  convex on $\R^n$ for fixed $\mu$ and concave on $\oM$ for fixed
  $\vlambda$ we may use the Min-Max Theorem [H-L] to obtain
   \ \ \begin{multline} \oE_\Xi^*(\vj)= \sup_{\mu\in \oM}\inf_{\vlambda\in \R^n}
   \left[ \int_{{\ttt}}\Xi d\mu + \oF(\mu, \vlambda) - \vlambda\cdot
   \vj\right] \\ =
   \inf_{\vlambda\in \R^n}\left\{ \sup_{\mu\in \oM}
   \left[ \int_{{\ttt}}\Xi d\mu + \oF(\mu,\vlambda)\right] - \vlambda\cdot
   \vj\right\}= -\sup_{\vlambda\in \R^n}\left\{ \vlambda\cdot\vj - \oE_\Xi(\vlambda)\right\}\end{multline}
   \ \ \ $\Box$
  \subsection{Proof of Proposition~\ref{cordelta2}}
  Let $\Xi(x_0)=\max_{{\ttt}} \Xi(x)$. Let
  $\mu_n\rightarrow\delta_{x_0}$. By Lemma\ref{fund}-(\ref{usc}, \ref{usefordelta}) and (\ref{oG}) we
  obtain the weak inequality for $\oE_\Xi(\vlambda)$. To obtain the
  strong inequality use, e.g., the uniform Lebesgue  measure
  $\mu=\dx$ on $\ttt$. Then the minimizer $\phi$ of $\oF$
  (\ref{Fbar2}) verifies $\nabla\cdot(\nabla\phi+\vlambda)=0$, that
  is, $\Delta\phi=0$ on $\ttt$ which implies $\nabla\phi\equiv 0$.
  Hence $\oF(\dx,\vlambda)= |\vlambda|^2/2$. We obtain the strong
  inequality if $|\vlambda|^2/2 +\int_{{\ttt}} \Xi > \max_{{\ttt}} \Xi$. The
  proof for $\oE^*$ is analogous.
  \ \ \ \ $\Box$

  \subsection{Proof of Lemma~\ref{muexists}}   The existence of a
  maximizer $\mu$ of (\ref{oG*}) follows from the
  lower-semi-continuity of $\oF^*$ (hence of $H_\Xi^*$) with respect to
  the $C^*$ topology of $\oM$, as claimed in
  Lemma~\ref{fund}-(\ref{usc}).
  \par
  Let  \be\label{11} H_\Xi( \mu, \vlambda):=  \int \Xi d\mu +  \oF(\mu, \vlambda) \ , \ \ \ H_\Xi^*( \mu, \vj):=
   \int \Xi d\mu -  \oF^*(\mu, \vj) \ . \ee Then, by
Definition~\ref{G}, $\oE_\Xi(\vlambda)= H_\Xi(\mu_0,\vlambda)\geq
H_\Xi(\mu,\vlambda)$ for any $\mu\in \oM$, and
$\oE_\Xi^*(\vj)=\sup_{\mu\in\oM}H_\Xi^*(\mu,\vj)$. By
Lemma~\ref{dualkj}
$$\oE_\Xi(\vlambda) - \oE_\Xi^*(\vj^{'}) \geq \vlambda\cdot \vj^{'}$$
holds for any $\vj^{'}\in\R^n$, and the equality above takes place
if and only if $\vj^{'}=\vj\in\partial_{\vlambda}\oE_\Xi(\cdot)$.
hence  \be\label{22} H_\Xi(\mu_0, \vlambda) - H_\Xi^*(\mu, \vj^{'})
\geq \vlambda\cdot \vj^{'}\ee holds for any $\mu\in\oM$ and any
$\vj^{'}\in\R^n$. The equality holds if and only if $\vj^{'}=\vj\in
\partial_{\vlambda}\oE_\Xi(\cdot)$ and $\mu$ which verifies the maximum of $H_\Xi^*(\cdot, \vj)$.
From (\ref{11}, \ref{22}) \be\label{33} \oF(\mu_0, \vlambda) +
\oF^*(\mu, \vj^{'}) +\int_{\mathbb{T}^n}\Xi( d\mu_0-d\mu) \geq
\vlambda\cdot \vj^{'} \ . \ee Let now $\mu=\mu_0$. Then (\ref{33})
implies \be\label{44} \oF(\mu_0, \vlambda) + \oF^*(\mu_0, \vj^{'})
\geq \vlambda\cdot \vj^{'} \ , \ee and, {\it if} \ there is an
equality in (\ref{44}) for some $\vj^{'}$, then $\vj^{'}\in
\partial_{\vlambda}\oE_\Xi(\cdot)$. However, we know, by definition
of $\oF^*$ as the Legendre transform of $\oF$ with respect to
$\vlambda$, that there is, indeed, an equality in (\ref{44})
provided $\vj^{'}\in\partial_{\vlambda}\oF(\mu_0,\cdot)$. This
verifies $\partial_{\vlambda}\oF(\mu_0,\cdot)\subset
\partial_{\vlambda}\oE_\Xi$.
  \subsection {Proof of
\rth{connect}}  Assume $\nu\in{\cal M}_L$ is a  Mather measure on
$\mathbb{T}^n\times\R^n$. Theorem 5.1.2 of [F] implies the existence
of a conjugate pair $\phi_+ \geq \phi_- \in Lip(\ttt)$ where the
domain $\phi_+=\phi_-$ contains the projected Mother set which, in
turn, contains the support of the projection $\mu$ of $\nu$ on
$\ttt$. In addition, Corollary~4.2.20 of [F] implies that either
functions satisfies
\be\label{hj1}\frac{1}{2}|\nabla\phi+\vlambda|^2+\Xi\leq E , \ee and
for some $E\in\R$, with an equality on the projected Mather set (in
particular, on the support of $\mu$).

   We show that $\mu$ is also a maximizer of
$\oE_\Xi(\vlambda)$ (\ref{oG}).
\par
Let  $\phi$ be either $\phi_+$ or $\phi_-$. Let $\eta_\eps\in
C^\infty(\mathbb{T}^n)$ non-negative mollifier function on
$\mathbb{T}^n$, supported in the ball $|x|<\eps$ and satisfying
$\int_{\mathbb{T}^n} \eta_\eps = 1$. Let $\phi^\eps:=
\phi*\eta_\eps$. Then $\phi^\eps\in C^\infty(\mathbb{T}^n)$ and
$\nabla\phi^\eps=\eta_\eps * \nabla\phi$.  The Jensen's inequality
implies that \be\label{jensen} \eta_\eps
* |\nabla\phi+\vlambda|^2 \geq |\nabla\phi^\eps+\vlambda|^2 \ , \ee
so, by (\ref{hj1}, \ref{jensen}))
$$ \frac{1}{2}|\nabla\phi^\eps+\vlambda|^2+\Xi\leq E +\Xi-\eta_\eps*\Xi \ . $$
Given $\delta>0$, there exists $\eps>0$ for which
$|\Xi-\eta_\eps*\Xi| < \delta$ on $\mathbb{T}^n$. Hence
$$ \frac{1}{2}|\nabla\phi^\eps+\vlambda|^2+\Xi\leq E +\delta \ . $$
So, for any $\tilde{\mu}\in\oM$:
$$ E+\delta \geq \int_{\mathbb{T}^n}\left(
\frac{1}{2}|\nabla\phi^\eps+\vlambda|^2+\Xi\right) d\tilde{\mu} \geq
\oF(\tilde{\mu}, \vlambda) + \int_{\mathbb{T}^n}\Xi d\tilde{\mu} =
\oE_\Xi(\vlambda) \ , $$ where the second inequality follows from
the definition (\ref{Fbar2}) of $\oF(\mu, \vlambda)$ and from
$\phi^\eps\in C^1(\mathbb{T}^n)$. Since $\delta>0$ can be
arbitrarily small it follows that $\oE_\Xi(\vlambda) \leq E$. Hence
$\mu=\pi_\#\nu$ is, indeed, a maximizer of $\oE_\Xi$ (\ref{oG}).
\par Finally, if $\vj=\vec{\alpha}(\nu)$ then
$\vj\in\partial_{\vlambda}\oE_\Xi$ and the last part of the Theorem
follows from Lemma~\ref{muexists}. \ \ \ \ \ $\Box$
\subsection{Proof of Proposition~\ref{thT}}
 The proof is analogous to this of \rth{th1}, utilizing
Proposition~\ref{propdual}. We only sketch the new definitions
involved, generalizing those given in the proof of \rth{th1} to the
time periodic case.
\par
Let
 ${\bf C}$ to be the space of all continuous functions on
$\ttt\times\R^n\times [0,T]$, equipped with the norm
$$\| q\|:=  \sup_{(x,\vp,t)\in \ttt\times
\R^n\times [0,T]}\left\{\frac{|q(x,\vp,t)|}{1+|\vp|} \right\}  \ .
$$ Define  \be\label{ZdefT} {\bf Z}:= \left\{ \phi_t+
\vp\cdot\nabla_x\phi  \ \ , \ \ \phi\in \CT \  \right\} \ee The dual
space ${\bf C}^*$ {\it contains} the set $\oM_T$ of $\oM-$valued
orbits on $[0,T]$ of bounded first moment. If $\hat{\nu}\in{\bf
C}^*$ is such an orbit then the duality relation is given by
$$ \langle\hat{\nu}, q\rangle:= \int_0^T\int_{
\ttt\times\R^n} qd\nu_{(t)}dt  \  \ , \  \forall q\in {\bf C} \ . $$
Given an orbit  $\hat{\mu}\in \oM_T$, define
\begin{multline} {\cal F}_{\hat{\mu}}(\hat{\nu},T) :=
 \int_{\ttt\times\R^n\times[0,T]} \left( |\vp|^2/2 -\vlambda\cdot\vp\right)d\nu_{(t)}dt  \\  \text{if}
  \ \ \hat{\nu}\in {\cal M}_T(\ttt\times\R^n\times [0,T])\cap {\bf C}^* \
    \text{and} \ \nu_{(t)}(\dx, \R^n) = \mu_{(t)}(\dx)\ \  \  \ a.e
 \\ {\cal F}_{\hat{\mu}}(\hat{\nu},T)= \infty \ \ \text{otherwise}  \ .
\end{multline}
The rest of the proof is  equivalent to this of \rth{th1}.
 \ \ \ \ \ \ $\Box$
\subsection{Proof of Proposition~\ref{maxoverline}}  First we note
that $\oF_T(\mu,\vlambda)\geq 0$ for any $\mu\in\oM$ and
$\vlambda\in\R^n$. Indeed, by \rth{thT}:
\be\label{estim}\oF_T(\mu,\vlambda)\geq \oF(\mu,\vlambda)\geq 0 \ .
\ee Let $\hat{\mu}^{(n)}$ be a maximizing sequence of
(\ref{FbarTmin}). By (\ref{defFT}) and (\ref{estim}) it follows that
there exists $C>0$ for which \be\label{qqq}\frac{1}{T}\int_0^T\int_{
\ttt\times\R^n} \left( \frac{|\vp|^2}{2}-\vp\cdot\vlambda\right)
d\nu^{(n)}_{(t)} dt\leq C\ee where $\hat{\nu}^{(n)}$ are the weak
solutions corresponding to $\hat{\mu}^{(n)}$.
\par
 Let
$$ \|\hat{\mu}\|_T^2 := -\inf_{\phi\in C^1(\ttt\times[0,T])}
\int_0^T\left( \phi_t +
\frac{1}{2}|\nabla_x\phi|^2\right)d\mu_{(t)}dt \ . $$ We recall form
[W] that \be\label{defromW} \|\hat{\mu}\|_T^2=
\frac{1}{2}\inf_{\hat{\nu}\in \widehat{\Lambda}_{\hat{\mu}}}
\int_0^T\int_{\ttt\times\R^n} |\vp|^2 d\nu_{(t)} dt \ . \ee
Moreover, Lemma~2.2  in [W] implies that the set $\left\{
\hat{\mu}\in \widehat{\Lambda}_{\hat{\mu}} \ ; \ \
\|\hat{\mu}\|_T<C\right\}$ is uniformly bounded in the $1/2-$Holder
norm with respect to the Wasserstein metric:
$$ W_1(\mu_1,\mu_2):= \sup_{\phi\in C^1(\ttt), |\nabla\phi|\leq 1} \int_{{\ttt}}\phi(x)(
d\mu_1-d\mu_2)= \inf_{\sigma\in\overline{\cal M}^{(2)}(\mu)
}\int_{{\ttt}}\int_{{\ttt}} |x-y|_T d\sigma \ , $$ (recall points
(i) and (vi) in List of Symbols). Since The $W_1$-Wasserstein metric
is a metrization of the weak topology of measures on compact
domains, it follows (see Corollary~2.1 in [W]) that the set $\left\{
\hat{\mu}\in \widehat{\Lambda}_{\hat{\mu}} \ ; \ \
\|\hat{\mu}\|_T<C\right\}$ is pre-compact in the topology of
$C\left([0,T]; C^*(\ttt)\right)$.
\par
Next, since $|\vp|^2/2-\vlambda\cdot\vp \geq |\vp|^2/4 -
|\vlambda|^2$ it follows from (\ref{qqq}) and (\ref{defromW}) that
$$ C\geq \frac{1}{T}
\inf_{\hat{\nu}\in\widehat{\Lambda}_{\hat{\mu}^{(n)}}}
\int_0^T\int_{{\ttt}} \left( \frac{|\vp|^2}{2}-\vlambda\cdot
\vp\right) d\nu_{(t)}dt \geq
\frac{1}{2T}\|\hat{\mu}^{(n)}\|_T^2-|\vlambda|^2 \ . $$ Hence, the
limit
$$\hat{\overline{\mu}}=\lim_{n\rightarrow\infty} \hat{\mu}^{(n)}\in
\widehat{\Lambda}_{\hat{\mu}}$$ exists in the weak topology of
$C^*([0,T]; C^*(\ttt))$, along a subsequence of the maximizing
sequence. Moreover, $\hat{\overline{\mu}}$ is a maximizer of
(\ref{FbarTmin}) by concavity of $\oF(\hat{\mu}, \vlambda,T)$ with
respect to $\hat{\mu}$.
 \ \ \ \ \ $\Box$
 \par\noindent
\subsection{Proof of Proposition~\ref{prop4.3}} The inequality
(\ref{obvineq}) follows directly from (\ref{FbarTmin}), upon
substituting the {\it constant} orbit $\hat{\mu}\equiv \mu_{(t)}
\equiv \mu$ for $t\in[0.T]$. To verify (\ref{oF=oG}) we use
(\ref{defFT}) to write \be\label{sumnumu} \sup_{\mu\in
\oM}H_{\Xi,T}(\mu,\vlambda)= \sup_{\hat{\nu}\in \widehat{\Lambda}_T}
\frac{1}{T}\int_0^T\int_{\ttt\times\R^n}\left\{
\Xi-\left[\frac{|\vp|^2}{2}-\vp\cdot\vlambda\right]
\right\}d\nu_{(t)}dt \ . \ee But, the RHS of (\ref{sumnumu}) is
unchanged if we replace $\hat{\nu}$ by $\nu:= T^{-1}\int_0^T
d\nu_{(t)}$ Moreover, $\nu\in \Lambda$ since Definition~\ref{def4.5}
is reduced to Definition~\ref{deflift} for $\theta(x,t) \rightarrow
\theta(x):= T^{-1}\int_0^T \theta(x,t)dt$. Hence (\ref{sumnumu}) is
reduced into
$$ \sup_{\mu\in \oM}H_{\Xi,T}(\mu,\vlambda)= \sup_{\nu\in
\Lambda} \int_{\ttt\times\R^n}\left\{
\Xi-\left[\frac{|\vp|^2}{2}-\vp\cdot\vlambda\right] \right\} d\nu \
,
$$ which yields (\ref{oF=oG}) via (\ref{defFT}). \ \ \ \ $\Box$
\subsection{Proof of Lemma~\ref{masstrans}}
The lower-semi-continuity of $D^T_{\vlambda}$ with respect to
$C^*(\oM)$ follows from the dual formulation
$$ D^T_{\vlambda}(\mu)=\sup_{\psi_1,\psi_2} \left( \int_{\mathbb{T}^n}
\psi_1 d\mu +\int_{\mathbb{T}^n} \psi_2 d\mu\right)$$ where the
supremum
above is taken over all pairs $\psi_1,\psi_2\in C(\mathbb{T}^n)$ verifying \\
$\psi_2(y)+\psi_1(x)\leq A^\Xi_{\vlambda}(x,y,T)$ for any $x,y\in
\mathbb{T}^n$. For details, see [V], ch. 1.
\par
The upper-semi-continuity follows directly from definition
(\ref{dvlambda1}). Indeed, let $\sigma_j\in \oM^{(2)}(\mu)$ verifies
(\ref{dvlambda1}) for $\mu_i\in\oM$ and $\mu_j\rightarrow\mu$ in
$C^*(\oM)$, then the sequence $\sigma_j$ is compact in the set
${\cal M}(\mathbb{T}^n\times \mathbb{T}^n)$ in the weak topology.
Let $\sigma$ be a limit of this sequence. Then
$\sigma\in\oM^{(2)}(\mu)$, and
$$ D^T_{\vlambda}(\mu)\leq \int_{\mathbb{T}^n} A^\Xi_{\vlambda}(x,y,T)
d\sigma =
\lim_{j\rightarrow\infty}\int_{\mathbb{T}^n}A^\Xi_{\vlambda}(x,y,T)
d\sigma_j = \lim_{j\rightarrow\infty}D^\Xi_{\vlambda}(\mu_j) \ . $$
The continuity of $D^\Xi_T$ with respect to $\Xi$ in the
$C^0(\mathbb{T}^n)$ topology is verified by the continuous
dependence of $A^\Xi_T$ on $\Xi$.  \ \ \ \ $\Box$
\subsection{Proof of Proposition~\ref{FtoD1}}
Definition~\ref{DefFT} of $\oE_T$ corresponds to  Definition 3.1 of
${\cal L}$ in [W].  In addition, Definition 3.3 [W] of ${\cal K}$
corresponds to (\ref{dvlambda}). Then, Proposition~\ref{DefFT}  is a
result of the identity ${\cal L}(\mu_1,\mu_2)= {\cal
K}(\mu_1,\mu_2)$, which follows from the Main Theorem of [W].
\par
In fact, the extended  Lagrangian ${\cal L}$ is defined, in [W], for
a Lagrangian $L=|\vp|^2/2 - \Xi(x)$, i.e. for $\vlambda=0$, but the
proof of the Main Theorem in [W] can be extended to $L=
|\vlambda-\vp|^2/2-\Xi$ in a direct way. \ \ \ \ \ $\Box$

\subsection{Proof of Proposition~\ref{FtoD2}} Since
$\oE_T(\mu,\vlambda, \Xi) \geq \int_{\tt}\Xi d\mu + \oF(\mu,
\vlambda)$ by Definition~\ref{DefFT}, it follows from
Proposition~\ref{FtoD1} that
$$  \frac{|\vlambda|^2}{2}-
D_{\vlambda}^T(\mu, \Xi) \geq \int_{\ttt}\Xi d\mu + \oF(\mu,
\vlambda) \ . $$ Thus, we only need to show that
$$ \lim_{T\rightarrow 0}  D_{\vlambda}^T(\mu,\Xi)\geq
\frac{|\vlambda|^2}{2} -\int_{\ttt}\Xi d\mu -\oF(\mu, \vlambda) \ .
$$ We may reduce to the case $\Xi\equiv 0$, hence we need to verify
\be\label{finresprop2} \lim_{T\rightarrow 0} \frac{|\vlambda|^2}{2}
- D_{\vlambda}^T(\mu)\leq
  \oF(\mu, \vlambda) \ . \ee Indeed, from
(\ref{Adef}, \ref{dvlambda}) we observe
$$ \lim_{T\rightarrow 0} D^T_{\vlambda}(\mu,\Xi) = -\int_{\ttt} \Xi d\mu + \lim_{T\rightarrow 0}
D^T_{\vlambda}(\mu) \ . $$
 If $\sigma_T$ verifies the minimum
in (\ref{dvlambda}), then \be\label{xyestim}D_{\vlambda}^T(\mu)=
\frac{1}{2T^2}\int_{\ttt}\int_{\ttt} \|\vy-\vx+\vlambda T\|^2
d\sigma_T(\vx,\vy)    \leq \frac{|\vlambda|^2}{2} \ . \ee Let
$$G:= \left\{ (x,y)\in \ttt\times\ttt \ ; \ \ \|y-x\|_\infty\geq
1/3\right\}  \ , $$ where $\|x\|_\infty$ is the metric on $\ttt$
defined as $\min_{\vz\in \mathbb{Z}}|\vx-\vz|_\infty$.  Then, for
sufficiently small $T$, $\| y-x+\vlambda T\|^2\geq 1/16$ for
$(x,y)\in G$ so from (\ref{xyestim}) \be\label{GTestim}
\int\int_{G}d\sigma_T(x,y)\leq 16T^2|\vlambda|^2  \ . \ee Let
$B^n\subset \R^n$ be the unite box $-1/2\leq x_i\leq 1/2$,
$i=1,\ldots n$.   Let $\phi\in C_0\left(\ttt\times T^{-1}
B^n/3\right)$. Extend $\phi$ to a function in $C_0\left(\ttt\times
T^{-1} B^n\right)$ by $\phi(x,\vp)=0$ if $\vp\in
T^{-1}B^n-T^{-1}B^n/3$. Further, extend $\phi$ into a function on
$\ttt\times \R^n$ as a $T^{-1}$ periodic function in $\vp$, that is,
$\phi$ is a function on $\ttt\times (\ttt/T)$. Set
$$ y=x+\vp T \ \ , \ \ \hat{\phi}(x,y):= \phi\left(x,
\frac{y-x}{T}\right)$$ Then $\hat{\phi}\in C(\ttt\times\ttt)$. Given
$\sigma_T$ which verifies the minimum in (\ref{dvlambda}), we define
a corresponding measure $\nu_T$ on $\ttt\times\R^n$, supported in
$\ttt\times (B^n/(3T))$, as follows: For any $\phi\in
C_0\left(\ttt\times  B^n/(3T)\right)$,
\be\label{nuTdef}\int_{\ttt}\int_{B^n/(3T)}\phi(x,\vp) d\nu_T(x,\vp)
= \int_{\ttt}\int_{\ttt}
 \hat{\phi}\left( x, y\right) d\sigma_T(x, y) \ . \ee By
 (\ref{GTestim}) we obtain
 \be\label{closeto1} 1\geq \int_{\ttt}\int_{\R^n} d \nu_T \geq 1- 16T^2
|\vlambda|^2 \ . \ee We now verify that \be\label{limnuT}
\lim_{T\rightarrow 0} \nu_T= \nu_0\in \Lambda_\mu\ee (see
Definition~\ref{deflift}-(2)). First, we show that the sequence of
measures on $\ttt\times\R^n$ is tight. For this, we use
(\ref{nuTdef}) with $\phi(x,\vp)=|\vp-\vlambda|^2 \cdot
1_{B^n/(3T)}(\vp)$ and (\ref{xyestim}) to obtain
$$ \frac{1}{2}\int_{\ttt}\int_{\R^n} |\vp-\vlambda|^2 d\nu_T < \frac{|\vlambda|^2}{2} \ . $$
This, and (\ref{closeto1}), imply that the limit (\ref{limnuT}) is a
probability measure on  $\ttt\times\R^n$, that is $\nu_0\in {\cal
M}(\ttt\times\R^n)$. Moreover, $\nu_0(dx, \R^n)= \mu(dx)$. To show
that $\nu_0\in\Lambda_\mu$, we proceed as follows: Let $q\in
C^\infty (\R^+)$ satisfies:
\begin{description}
\item{i)} \  $q\in C^\infty(\R^+)$.
\item{ii)} \ $q(s)=1$ for $0\leq s\leq 1/3$.
\item{iii)} \ $q(s)=0$ for $ 2/5\leq s\leq 1/2$.
\item{iv)} \ $q(1/2-s)=q(1/2+s)$ for any $s\in[0,1/2]$.
\item{ v)} \ $q(s)\leq 1$ for $s\in [0,1]$.
\item{vi)} \ $q(s+1)=q(s)$ for all $s\in\R^+$.
\end{description}
 Let
$$ Q(\vp)= \Pi_1^N q(|p_i|) \ \text{for} \ \vp=(p_1,
\ldots p_n)\in\R^n \ \ , \ \ \ Q_T(\vp)= Q(T\vp) \ .
$$
 Given $\psi\in C^1(\ttt)$, set $\phi^{(1)}_T(x,p)= \psi(x+pT)
Q_T(\vp)$ and $\phi^{(2)}_T(x,p)= \psi(x) Q_T(\vp)$ if $\vp\in
B^n/(3T)$, $\phi_T=0$ otherwise. Then, by (\ref{nuTdef})
\begin{multline}\label{phi1}
 \int_{\ttt}\int_{\R^n} \phi^{(1)}_T(x,p)
d\nu_T(x,p)=\int_{\ttt}\int_{\ttt} Q(\vy-\vx) \psi(y)
d\sigma_T(\vx,\vy) = \\ \int_{\ttt}\int_{\ttt} \psi(y)
d\sigma_T(\vx,\vy) + \int_{\ttt}\int_{\ttt}( Q(\vy-\vx)-1) \psi(y)
d\sigma_T(\vx,\vy)  \\
 =   \int_{\ttt} \psi(y) d\mu(\vy) +
 \int_{\ttt}\int_{\ttt}( Q(\vy-\vx)-1) \psi(y)
d\sigma_T(\vx,\vy) \ .
 \end{multline}
 The same argument applies also to $\phi^{(2)}$ and yields
\begin{multline}\label{phi2}
\int_{\ttt}\int_{\R^n} \phi^{(2)}(x,p) d\nu_T(\vx,\vp)
 =   \int_{\ttt} \psi(x) d\mu(\vx) +
 \int_{\ttt}\int_{\ttt}( Q(\vy-\vx)-1) \psi(x)
d\sigma_T(\vx,\vy) \ .
 \end{multline}
 However, $Q(\vy-\vx)-1=0$ on the set $G$ so, by
 (\ref{GTestim})
\be\label{correct} \int_{\ttt}\int_{\ttt}( Q(\vy-\vx)-1) \psi(x)
d\sigma_T(\vx,\vy)=  \int_{\ttt}\int_{\ttt}( Q(\vy-\vx)-1) \psi(y)
d\sigma_T(\vx,\vy)\leq 16T^2|\vlambda|^2|\psi|_\infty \ . \ee
Subtract (\ref{phi2}) from (\ref{phi1}), divide  by $T$ and let
$T\rightarrow 0$ and use (\ref{correct}) to obtain \begin{multline}
0= \lim_{T\rightarrow 0} \int_{\ttt}\int_{\R^n}
\frac{\phi^{(2)}(\vx,\vp)-\phi^{(1)}(\vx,\vp)}{T} d\nu_T(\vx,\vp) =
\\ \lim_{T\rightarrow 0} \int_{\ttt}\int_{\R^n}
Q_T(\vp)\frac{\psi(\vx+T\vp)-\psi(\vx)}{T} d\nu_T(\vx,\vp)=
 \int_{\ttt}\int_{\R^n}
\nabla\psi\cdot\vp d\nu_0(\vx,\vp)
 \ ,
\end{multline}
which implies for any $\psi\in C^1(\ttt)$, hence
$\nu_0\in\Lambda_\mu$ as claimed.
\par
Let now consider $\int_{\ttt}Q_T(\vp)|\vp-\vlambda|^2d\nu_T \ . $ By
(\ref{nuTdef})
\be\label{one}\frac{1}{2}\int_{\ttt}\int_{B^n/(3T)}Q_T(\vp)|\vp-\vlambda|^2
d\nu_T(x,\vp) = \frac{1}{2}T^{-2}\int_{\ttt}\int_{\ttt} Q(\vx-\vy)
\| \vx-\vy-T\vlambda\|^2d\sigma_T(x, y)  \ . \ee Since $\sigma_T$ is
a minimizer of (\ref{DTk0}) and $Q_T\leq 1$, it follows that
\be\label{two} \frac{1}{2}T^{-2}\int_{\ttt}\int_{\ttt} Q(\vx-\vy) \|
\vx-\vy-T\vlambda\|^2d\sigma_T(x, y)\leq D_{\vlambda}^T(\mu) \ . \ee
On the other hand, $\nu_T$ is supported on $B^n/(3T)$ by definition,
so
\be\label{three}\frac{1}{2}\int_{\ttt}\int_{B^n/(3T)}Q_T(\vp)|\vp-\vlambda|^2
d\nu_T(x,\vp)=\frac{1}{2} \int_{\ttt}\int_{\R^n}|\vp-\vlambda|^2
d\nu_T(x,\vp) \ . \ee On the other hand \be\label{four}
\lim_{T\rightarrow 0} \int_{\ttt}\int_{\R^n}|\vp-\vlambda|^2
d\nu_T(x,\vp)\geq\int_{\ttt}\int_{\R^n}|\vp-\vlambda|^2
d\nu_0(x,\vp) \ . \ee From(\ref{one}-\ref{four}) we obtain
$$ \frac{1}{2}\int_{\ttt}\int_{\R^n}|\vp-\vlambda|^2
d\nu_0(x,\vp)\leq \lim_{T\rightarrow 0} D_{\vlambda}^T(\mu) \ , $$
hence \be\label{th1fin} \lim_{T\rightarrow 0}  \left[
\frac{|\vlambda|^2}{2} - D_{\vlambda}^T(\mu)\right] \leq
-\int_{\ttt}\int_{\R^n}\left( \frac{|\vp|^2}{2}-\vlambda\cdot
\vp\right) d\nu_0(x,\vp)  \ . \ee Since we already proved that
$\nu_0\in \Lambda_\mu$, then \rth{th1} and (\ref{th1fin}) verify
(\ref{finresprop2}). \ \ \ \ \ $\Box$ \vskip .2in\noindent
\subsection{Proof of Theorem~\ref{DisH}}
To prove (\ref{athdish}) we use (\ref{oF=oG}) and (\ref{BBthC})
together with \be\label{Dish1}\lim_{j\rightarrow\infty}
D^T_{\vlambda}(j, \Xi)=\min_{\mu\in\oM} D^T_{\vlambda}(\mu, \Xi) \ .
\ee To establish (\ref{Dish1}) note, first, that $D^T_{\vlambda}(j,
\Xi)\geq \inf_{\mu\in\oM} D^T_{\vlambda}(\mu, \Xi)$ for any $j$ \ by
definition, so it is enough to establish the inequality
$$ \limsup_{j\rightarrow\infty}
D^T_{\vlambda}(j, \Xi)\leq \inf_{\mu\in\oM} D^T_{\vlambda}(\mu, \Xi)
\ . $$ Let now $\{ \mu_j\}$ be a sequence of empirical measures,
where, for each $j$, $\mu_j$ contains exactly $j$ atoms, and so that
$\mu=\lim_{j\rightarrow\infty} \mu_j$ in $C^*$. Then
$D^T_{\vlambda}(\mu_j, \Xi)\geq D^T_{\vlambda}(j, \Xi)$ by
definition, while $\lim_{j\rightarrow\infty}
D^T_{\vlambda}(\mu_j,\Xi)= D^T_{\vlambda}(\mu,\Xi)$ by
Lemma~\ref{masstrans}.

\begin{center}{\bf References}
\end{center}
\begin{description}
\item{[A]} S. Aubry: {\it The twist map, the extended
Frenkel-Kontrovna model and the devil's staircase}, Physica 7{\bf
D}, 240-258, 1983.
\item{[AS]} E. A. Arriola; J. Soler:  {\it A Variational Approach to the Schr$\ddot{\text{o}}$dinger-Poisson
System: Asymptotic Behaviour, Breathers, and Stability},  J. Stat.
Phys. {\bf 103}, no. 5-6, 1069-1106, 2001.
\item{[BB]} P. Bernard and B. Buffoni , {\it Optimal mass
transportation and Mather Theory}, Preprint, \ \
http://arxiv.org/abs/math.DS/0412299
\item{[E]} L. C. Evans, {\it  A survey of partial differential equations methods in weak KAM theory},
 Comm. Pure Appl. Math. {\bf 57}  4, 445-480, 2004
\item{[EG1]} L.C.Evans and  D. Gomes, {\it
Effective Hamiltonians and averaging for Hamiltonian dynamics. I.},
Arch. Ration. Mech. Anal. {\bf 157}  1, 1-33, 2001.
\item{[EG2]} L.C.Evans and  D. Gomes, {\it Effective Hamiltonians and averaging for
Hamiltonian dynamics. II}. Arch. Ration. Mech. Anal. {\bf 161} ,  4,
271-305, 2002
\item{[F]} \ A. Fathi: {\it The weak Kam Theorem in Lagrangian
Dynamics}, Cambridge University Press, Cambridge Studies in Advanced
Mathematics Series,  Volume {\bf 88}, 2003
\item{[G]} L Granieri: {\it On action minimizing measures for the
Monge-Kantorovich problem} Preprint,  (july 2004).
\item{[GO]}
D.A. Gomes and A.M.  Oberman: {\it  Computing the Effective
Hamiltonian Using a Variational Approach},  SIAM J. Control Optim.
{\bf 43}, 3, 792-812, 2004
 \item{[GCJ]} \ M. Govin, C.
Chandre and H.R. Jauslin: {\it
Kolmogorov-Arnold-Moser-Renormalization-Group analysis of stability
in Hamiltonian flows}, Phys. Rev. Lett. {\bf 79}, 3881-3884, 1997
\item{[GT]} D. Gilbarg, N. S. Trudinger.: {\it Elliptic partial differential
equations of second order},
 Berlin ; New York : Springer-Verlag, [1977].
\item{[H]} \ G.A.
Hedlund, {\it Geodesics on a 3 dimensional Riemannian manifolds with
periodic coefficients}, Ann Math. {\bf 33}, 719-739, 1932
\item{[H-L]} J.B. Hiriart-Urruty and C. Lemar¶echal: {\it Convex
Analysis and Minimization Algorithms II}, volume {\bf 306} of
Grundlehren der Mathematischen Wissenschaften, chapter 10.
Springer-Verlag, 1993.
\item{[IZL]}
R. Illner, P. F. Zweifel, and H. Lange: {\it Global existence,
uniqueness and asymptotic behaviour of solutions of the
Wigner–Poisson and Schr$\ddot{\text{o}}$dinger–Poisson systems},
Math. Meth. Appl. Sci. {\bf 17}, 349-376, 1994.
\item{[K]} J. B. Keller: {\it Semiclassical mechanics}, SIAM Review, {\bf 27},
4, 485-504, 1985
\item{[LSG]} D. P. Luigi, G.M. Stella and L. Granieri: {\it Minimal measures, one-dimensional currents and the
Monge-Kantorovich problem,} Calc. Var. Partial Differential
Equations {\bf 27}, 1, 1-23, 2006
\item{[Man]} R. Ma$\tilde{n}\grave{e}$: {\it On the minimizing
measures of Lagrangian dynamical systems}, Nonlinearity {\bf 5},
623-638, 1992 \item{[Mat]} J.N. Mather: {\it Existence of
quasi-periodic orbits for twist homeomorphisms on the annulus},
Topology {\bf 21}, 457-467, 1982 \item{[Mat1]} J.N. Mather: {\it
Minimal measures}, Comment. Math. Helv. , {\bf 64}, 375-394, 1989
\item{[Mo]} J. Moser: {\it Monotone twist mappings and the
calculus of variations}, Ergod. Th. \& Dynam. Sys., {\bf 6},
401-413, 1986 \item{[RW]} J. Rubinstein and G. Wolansky, {\it
Eikonal functions: Old and new}, in A Celebration of Mathematical
Modeling: The Joseph B. Keller Anniversary Volume, Kluwer, D.
Givoli, M.J. Grote,  G. Papanicolaou, eds., 2004.
\item{[S]} K.F. Siburg: {\it The Principle of
Least Action in geometry and Dynamics} , Lecture Notes in
Mathematics {\bf 1844}, Springer,2004.
\item{[V]} C. Villani: {\it Topics in Optimal Transportation},
Graduate studies in Math, {\bf 58}, AMS, 2003
\item{[W]} G. Wolansky: {\it Optimal Transportation in the presence of a
prescribed pressure field} , Preprint \ \ arXiv:math-ph/0306070 v5
12 Jan 2006
\end{description}

\end{document}